\documentclass[preprint, 12p]{imanum}

\makeatletter
\def\ps@plain{%
	\let\@oddhead\@empty
	\let\@evenhead\@empty
	\def\@oddfoot{}%
	\let\@evenfoot\@oddfoot}
\def\ps@myheadings{%
	\let\@oddhead\@empty
	\let\@evenhead\@empty
	\def\@oddfoot{}%
	\let\@evenfoot\@oddfoot}
\def\@maketitle{%
	\newpage\setcounter{footnote}{1}
	\null
	\vspace*{12pt}%
	\parbox{132mm}{
		\begin{center}%
			\let \footnote \thanks
			{\fontsize{12}{15}\bf\selectfont \@title \par}%
			\vskip 1.2em%
			{
				{\fontsize{10}{12}\it\selectfont \@author\\[6pt]
					\phantom{\rm[\@received\@revised\@accepted]}\\[9.5pt]
				}
			}
			\vskip 1em%
		\end{center}
	}
	\vspace*{-.75pc}}
\makeatother

\usepackage[T1]{fontenc}
\usepackage{graphicx}
\usepackage{lipsum}
\usepackage{amsfonts}
\usepackage{graphicx}
\usepackage{amsfonts}
\usepackage{mathtools}
\usepackage{hhline}
\usepackage{array}
\usepackage{makecell}
\usepackage[caption=false]{subfig}
\usepackage{bbm}
\usepackage{placeins}
\usepackage{cleveref}

\renewcommand{\cite}[2][]{\citep[#1]{#2}}

\usepackage{tikz}
\usepackage{pgf}
\usepackage{pgfplots}
\usepgfplotslibrary{groupplots}
\pgfplotsset{compat=newest}
\usetikzlibrary{spy}
\usetikzlibrary{shapes}
\usetikzlibrary{backgrounds}
\usetikzlibrary{shadows}
\usetikzlibrary{matrix}
\usepackage{shellesc}
\pgfplotsset{compat=newest}
\usetikzlibrary{external}
\usetikzlibrary{positioning}
\usetikzlibrary{plotmarks}
\usepgfplotslibrary{fillbetween}
\usetikzlibrary{decorations.pathmorphing}
\usetikzlibrary{matrix}
\usetikzlibrary{arrows} 
\usetikzlibrary{calc}
\usetikzlibrary{shapes}

\tikzstyle{block} = [draw,rectangle,thick,minimum height=2em,minimum width=2em]
\tikzstyle{sum} = [draw,circle,inner sep=0mm,minimum size=2mm]
\tikzstyle{connector} = [->,thick]
\tikzstyle{line} = [thick]
\tikzstyle{branch} = [circle,inner sep=0pt,minimum size=1mm,fill=black,draw=black]
\tikzstyle{guide} = []
\tikzstyle{snakeline} = [connector, decorate, decoration={pre length=0.2cm,
	post length=0.2cm, snake, amplitude=.4mm,
	segment length=2mm},thick, magenta, ->]

\allowdisplaybreaks

\newcommand{\x}{\vec{x}}

\newcommand{\op}[1]{\operatorname*{#1}}
\newcommand{\eps}{\varepsilon}
\newcommand{\norm}[1]{\left\lVert #1 \right\rVert}
\newcommand{\abs}[1]{\left| #1 \right|}
\newcommand{\dt}{h} 

\newcommand{\sca}[2]{\langle #1, \; #2 \rangle}

\renewcommand{\exp}[1]{e^{#1}}

\renewcommand{\d}{\operatorname{d}\!}

\newcommand{\Ind}{\mathbbm{1}}
\renewcommand{\O}{\mathcal{O}}
\newcommand{\conv}{*}

\renewcommand{\vec}[1]{\pmb{#1}}

\newcommand{\drv}[2]{\partial_{#1} #2}

\newcommand{\grad}{\nabla}

\newcommand{\Ea}{\op{E}\!} 

\newcommand{\LaSymb}{\mathfrak{L} }
\newcommand{\La}[1]{\LaSymb\left\lbrace #1 \right\rbrace }

\newcommand{\C}{\mathbb{C}} 

\newcommand{\HS}{\mathcal{H}} 
\renewcommand{\L}[1]{\mathrm{L}^{#1}} 
\renewcommand{\H}[1]{\mathrm{H}^{#1}} 
\newcommand{\Bochner}{\mathcal{B}} 
\newcommand{\BochnerDual}{\hat{\mathcal{B}}} 

\newcommand{\fenics}{FEniCS}

\newcommand{\FrDerRL}[1]{\partial^{#1}_t}
\DeclareMathAlphabet{\mathscr}{OMS}{cm}{m}{n}
\newcommand{\FrInt}[1]{\mathscr{I}^{#1}}
\newcommand{\FrKer}[1]{K_{#1}}
\newcommand{\ExpKer}[1]{\tilde{K}_{#1}}

\newcommand{\scaH}[2]{\sca{#1}{#2}_{\HS}}
\newcommand{\dualH}[2]{\sca{#1}{#2}_{\HS^\prime,\HS}}
\newcommand{\normH}[1]{\norm{#1}_{\HS}}
\newcommand{\normHdual}[1]{\norm{#1}_{\HS^\prime}}

\newcommand{\ErrRel}{\mathcal{E}_r}

\newcommand{\normB}[1]{\norm{#1}_{\Bochner}}
\newcommand{\normT}[1]{\normB{#1}}
\newcommand{\normBHdual}[1]{\norm{#1}_{\BochnerDual}} 
\newcommand{\normBA}[2][\alpha]{\norm{#2}_{\Bochner_{#1}}} 

\definecolor{color0}{rgb}{0.7843, 0.7843, 0.7843}
\definecolor{color1}{rgb}{0, 0.4470, 0.7410}
\definecolor{color2}{rgb}{0.8500, 0.3250, 0.0980}
\definecolor{color3}{rgb}{0.9290, 0.6940, 0.1250}
\definecolor{color4}{rgb}{0.7060, 0.3840, 0.7650}
\definecolor{color5}{rgb}{0.4660, 0.6740, 0.1880}
\definecolor{color6}{rgb}{0.3010, 0.7450, 0.9330}
\definecolor{color7}{rgb}{0.6350, 0.0780, 0.1840}
\definecolor{color8}{rgb}{0.0, 0.4078, 0.3412}

\pgfplotsset{
	legend image with text/.style={
		legend image code/.code={%
			\node[anchor=center] at (0.3cm,0cm) {#1};
		}
	},
}

\graphicspath{{./figures/}}

\title{Solving time-fractional differential equations via rational approximation}
\shorttitle{Solving time-FDE via rational approximation}

\author{%
	{\sc
		Ustim~Khristenko\thanks{Corresponding author. Email: khristen@ma.tum.de}
		and
		Barbara~Wohlmuth\thanks{Email: wohlmuth@ma.tum.de}
	} \\[2pt]
	Lehrstuhl f\"ur Numerische Mathematik, 
	Technische Universit\"at M\"unchen
}
\shortauthorlist{U.~Khristenko and B.~Wohlmuth}

\begin{document}
	
	\maketitle
	
	\begin{abstract}
		{Fractional differential equations (FDEs) describe subdiffusion behavior of dynamical systems.
		Its non-local structure requires taking into account the whole evolution history during the time integration, which then possibly causes additional memory use to store the history, growing in time.
		An alternative to a quadrature for the history integral is to approximate the fractional kernel with the sum of exponentials, which is equivalent to considering the FDE solution as a sum of solutions to a system of ODEs.
		One possibility to construct this system is to approximate the Laplace spectrum of the fractional kernel with a rational function.
		In this paper, we use the adaptive Antoulas--Anderson (AAA) algorithm for the rational approximation of the kernel spectrum which yields only a small number of real valued poles.
		We propose a numerical scheme based on this idea and study its stability and convergence properties.
		In addition, we apply the algorithm to a time-fractional Cahn-Hilliard problem.}
		{time-fractional differential equations, rational approximation, AAA algorithm.}
	\end{abstract}
	

	\section{Introduction}
	\label{sec:Intro}

Fractional differential equations have become more common components of models for complex physical systems in recent years as they often provide more realistic characterizations of certain physical phenomena than traditional differential operators.
In particular, time-fractional differential equations can be found in applications such as, e.g., modeling tumor growth and certain models in visco-elasticity.
In general, analytical solutions of such problems are not available.
Thus, the development and study of computational methods and algorithms for FDEs is of high interest and is the object of an increasing number of research works.
The numerical solution of FDEs is much more complicated and more computationally expensive than for the classical integer-order problem.
In fact, fractional ordinary differential equations (FODEs) can be regarded as integral equations involving convolution with singular kernels.
Thus, one of the related challenges is the non-local structure of the operators, which takes into account the whole evolution history during the time integration.
This feature causes additional memory to store the history, which grows in time.

The most common strategy is directly based on quadrature schemes for the convolution integral.
The classical one is the so-called L1 scheme based on a finite difference formula; see~\cite{oldham1974fractional,jin2016analysis}.
An important class constitutes the fractional linear multi-step methods (FLMM).
Pioneering work in this direction has been done by~\cite{lubich1983stability,lubich1986discretized,lubich1988convolution}.
The FLMM include methods of Adams--Moulton/Bashforth type; see, e.g.,~\cite{diethelm2002predictor,diethelm2004detailed,zayernouri2016fractional,zhou2020implicit}.
Though these methods are conceptually simple and have a high convergence order, they may experience difficulties for certain values of the fractional power; see, e.g.,~\cite{diethelm2006pitfalls}.
Moreover, the quadrature-based approach is mostly affected by the curse of non-locality.
Precisely, the integration of $N$ time steps has algorithmic complexity of order~$\O(N^2)$ and requires $\O(N)$ solutions to store.
The non-locality problem can be tackled with various memory-saving techniques including short memory principle~\cite{podlubny1998fractional,deng2007short}, logarithmic grids~\cite{ford2001numerical,diethelm2006efficient} and parallel computations~\cite{diethelm2011efficient} for the history integral.

Another strategy for the numerical solution of FODEs is based on the approximation of the integral kernel.
The first steps in this direction have been done in~\cite{lubich2002fast,schadle2006fast,lopez2008adaptive}.
In the so-called kernel compression technique, the kernel is approximated with a sum-of-exponentials~\cite{beylkin2005approximation,mclean2018exponential}, leading to a family of ODEs, which can be usually solved in parallel.
Therefore, this approach can be seen as a decomposition of the FODE solution into the sum of different modes, each governed by a corresponding local evolution law.
This phenomenon clearly illustrates the nature of FODEs, where the computational complexity can be interpreted in terms of a hidden extra-dimension, which is also reflected in the methods for fractional partial differential equations; compare, e.g., with \cite{caffarelli2007extension,banjai2019tensor,bonito2015numerical,vabishchevich2015numerically,harizanov2018positive}.
In this approach, the history term is approximated with a linear combination of $m=\O(\log N)$ auxiliary modes. 
Therefore, it also requires an additional memory storage of the size~$\O(\log N)$ and the computational complexity~$\O(N\,\log N)$.
In many works, the sum-of-exponentials is obtained using a quadrature for the integral representation of the kernel; see~\cite{li2010fast,mclean2006time,jiang2017fast,zeng2018stable,baffet2019gauss,banjai2019efficient}.
An alternative approach for the approximation of the kernel consists in polynomial (multi-pole) interpolation of its spectrum~\cite{baffet2017kernel}.
Our approach belongs to this class.
For the detailed overview of the existing numerical methods for FODE, we refer in particular to~\cite{baleanu2012fractional,diethelm2008investigation,diethelm2020good}.

In this work, we propose a new approach based on the approximation of the Laplace spectrum of the convolution kernel with a rational function.
We use the adaptive Antoulas--Anderson (AAA) algorithm~\cite{nakatsukasa2018aaa}.
The AAA algorithm was first applied to the solution of fractional diffusion problems in~\cite{hofreither2020unified}.
In the fractional ODE framework, it leads to a multi-pole approximation of the spectrum with real non-negative poles, which transforms to the sum-of-exponentials kernel with an additional singular term.
Discretizing the system of ODEs obtained for the modes, we propose new numerical schemes of the Implicit Euler and Crank-Nicolson types.
We also propose a stabilized version of the latter, based on an exponential integrator, which is able to avoid spurious oscillations typical for such schemes.
Moreover, we do not discretize explicitly the local integration term, which is obtained from the rational approximation.
Nevertheless, the coefficients of the local term in the numerical schemes naturally reproduce the fractional Adams--Moulton coefficients arising in linear multi-step methods.

In our method, though the number of modes $m$ grows as $\log N$ like in other methods, the value of $m$ is significantly less than the number of auxiliary variables (quadrature points) in many other kernel compression methods.
In particular, $m < 25$ is sufficient in our experience to achieve excellent accuracy even in complex non-linear PDE test cases.
This can play a decisive role in industrial applications using fine spatial discretizations in 3D, where the memory cost of each additional stored solution vector is tremendously high.

The principal novelty of our approach is that the sum-of-exponentials approximation can be constructed by fitting the Laplace spectrum of the fractional kernel along the real line only, while leading to a small number of modes.
Besides, the local term is constructed directly from the rational approximation of the spectrum. 

The paper is structured as follows.
In~\Cref{sec:Preliminaries}, we bring in the necessary definitions and formulations and provide some preliminary results.
In~\Cref{sec:RA}, we discuss the rational approximation of the spectrum of the fractional kernel.
In~\Cref{sec:EA}, we analyze the associated approximation error.
In~\Cref{sec:Scheme}, we suggest numerical schemes based on the discretization of a modal ODE system and discuss their accuracy and stability.
Finally, in~\Cref{sec:Examples}, we illustrate the performance of the newly introduced schemes numerically.
In particular, in~\Cref{sec:Linear}, we consider a linear time-fractional heat equation in 1D with known analytical solution to show the convergence rate.
In the second numerical example in~\Cref{sec:CahnHilliard}, the proposed scheme is applied to a non-linear time-fractional Cahn-Hilliard equation in~2D.

	\section{Preliminaries}
	\label{sec:Preliminaries}


Let us first introduce some basic definitions of the fractional derivative and the fractional integral.
For a more detailed introduction to the theory of fractional differential equations and fractional calculus, the reader is referred to~\cite{baleanu2012fractional,bajlekova2001fractional,diethelm2010analysis,kilbas2006theory,mainardi2010fractional,podlubny1998fractional,miller1993introduction,samko1993fractional}.
The Riemann--Liouville fractional integral is defined for $\alpha\in(0,\,1]$ as
\begin{equation}\label{eq:def:frac_integral}
	\FrInt{\alpha}u(t) = \FrKer{\alpha}*u(t) = \int_0^t \FrKer{\alpha}(t-s) u(s)\d s,
\end{equation}
where the kernel is defined by
\begin{equation}\label{key}
	\FrKer{\alpha}(t)=t^{\alpha-1}/\Gamma(\alpha),
\end{equation}
and $\Gamma(\alpha)$ is the gamma function.
Then, the Riemann--Liouville fractional derivative is given by
\begin{equation}\label{eq:FracDer_RL}
	\FrDerRL{\alpha} u(t)  = \partial_t \left[\FrInt{1-\alpha}u(t)\right].
\end{equation}
 If $u(t)$ is sufficiently smooth, then we have 
\begin{equation}\label{Def:frac_Caputo}
\begin{aligned}
\FrDerRL{\alpha}\left[u(t)-u(0)\right]=\FrInt{1-\alpha}\left[\partial_t u(t)\right].
\end{aligned}
\end{equation}
The right-hand side is the classical Caputo fractional derivative. The formulation on the left hand side, which expresses the Caputo fractional derivative in terms of Riemann--Liouville fractional derivative, has the advantage that it requires less regularity of $u(t)$ than the classical definition.


Let $\HS$ be a Hilbert space with inner product~$\sca{\cdot}{\cdot}$ and associated norm~$\normH{\cdot}$.
Let us consider the associated Bochner space~$\Bochner = \L{2}([0,T]; \HS)$ with the norm defined by $\normT{u}^2 = \int_{0}^{T}\normH{u(t)}^2\d t$.
In the current work, we focus on the solution $u\in\Bochner$ of the following non-linear fractional Cauchy problem:
\begin{align}\label{eq:NLequation}
\drv{t}{^\alpha u}(t) &= F[u](t), \\
u(0) &= u_0,
\end{align}
where $\alpha\in(0,\,1]$, and $F: \Bochner \to \BochnerDual$, with $\BochnerDual = \L{2}([0,T]; \HS^\prime)$, is a continuous possibly non-linear operator, such that $F[u](t) = F(t,u(t))$.
The solution of Equation~\eqref{eq:NLequation} can be formally written in terms of~\eqref{eq:def:frac_integral} as
\begin{equation}\label{eq:def:solution}
u(t) = u_0 + \FrKer{\alpha}*F[u].
\end{equation}
Let us denote by $\LaSymb$ the Laplace transform operator.
Then, for any function $f\in\Bochner$, we denote its Laplace transform by $\hat{f}(s) := \La{f}(s) = \int_{0}^{T}f(t)\exp{-st}\d t$, $s\in\C$.
In particular, the Laplace transform of~\eqref{eq:def:solution} is given by
\begin{equation}
\La{u-u_0}(s) = s^{-\alpha}\La{F[u]}(s).
\end{equation}

We are interested in the construction of an approximation of~$u(t)$ that we introduce through
\begin{equation}
\tilde{u}(t) = u_0 + \ExpKer{\alpha}*F[\tilde{u}],
\label{eq:def:approximation}
\end{equation}
where $\tilde{K}_{\alpha}(t)$ is an approximation of the kernel~$K_{\alpha}(t)$.


	\section{Rational approximation of the kernel}
	\label{sec:RA}

We want to construct a kernel~$\tilde{K}_{\alpha}(t)$ such that~\eqref{eq:def:approximation} yields a good approximation to
the solution~\eqref{eq:def:solution} of~\eqref{eq:NLequation}, which can be numerically found at a lower computational cost.
Let us consider the Laplace transforms of the kernels $K_{\alpha}(t)$ and $\tilde{K}_{\alpha}(t)$.
We look for~$\hat{\tilde{K}}_{\alpha}(z)$, $z\in\C$, in the form of a rational function, more precisely as a ratio of two polynomials.
That is, we want to construct $\hat{\tilde{K}}_{\alpha}(z)$ as a rational approximation of $\hat{K}_{\alpha}(z) = z^{-\alpha}$.
Since $\alpha\in(0,1]$, we select the numerator and the denominator of~$\hat{\tilde{K}}_{\alpha}(z)$ as polynomials of the same degree~$m$.
Under the assumption that the polynomials have only simple roots, the partial fractions decomposition of $\hat{\tilde{K}}_{\alpha}(z)$ is given by
\begin{equation}\label{eq:La_K_tilde}
	\hat{\tilde{K}}_{\alpha}(z) = \sum_{k=1}^{m} \frac{w_k}{z + \lambda_k} + w_{\infty},
\end{equation}
with $w_k\ge 0$, $\lambda_k \ge 0$.
Then, the kernel~$\tilde{K}_{\alpha}(t)$ is given by a typical sum-of-exponentials and a singular term:
\begin{equation}\label{eq:K_tilde}
\tilde{K}_{\alpha}(t) = \sum_{k=1}^{m} w_k\exp{-\lambda_k t} + w_{\infty}\delta(t),
\end{equation}
where $\delta(t)$ denotes the Dirac $\delta$-distribution.
Hence, the approximation~$\tilde{u}(t)$, defined by~\eqref{eq:def:approximation}, reads as
\begin{equation}\label{eq:RA_expansion}
\tilde{u}(t) = u_0 + \sum_{k=1}^{m} u_k(t) + u_{\infty}(t),
\end{equation}
where the modes~$u_k(t)$ are given by
\begin{align}
u_k(t) &= w_k\int_{0}^{t}\exp{-\lambda_k(t-s)}\,F[\tilde{u}](s)\d s, 
\qquad k=1,\ldots,m,\\
u_{\infty}(t) &= w_\infty F[\tilde{u}](t).
\label{eq:RA_modes}
\end{align}
That is, the modes~$u_k(t)$, $k=1,\ldots,m$, satisfy the following ordinary differential equation:
\begin{align}
\drv{t}{u_k}(t) + \lambda_k\,u_k(t) - w_k F[\tilde{u}] &= 0, \label{eq:sys_modes}
\\
u_k(0) &= 0.
\end{align}
Remark that the uni-modal ($m=1$, $w_\infty=0$) setting with $w_1=1$ and $\lambda_1=0$ yields the trivial case~$\alpha=1$.
On the other hand, the situation with only the "infinity" mode $w_{\infty}=1$, $m=0$, corresponds to the case~$\alpha=0$.
For the existence of $\tilde{u}$, which solves the Volterra integral equation~\eqref{eq:RA_expansion}-\eqref{eq:RA_modes}, we refer to~\cite{gripenberg1990volterra}.

Note that the expansion~\eqref{eq:RA_expansion}-\eqref{eq:RA_modes} is similar to the representations in~\cite{li2010fast,Yuan2002,diethelm2008investigation}; however, we obtain the weights~$w_k$ and the exponents~$\lambda_k$ from the rational approximation of the kernel spectrum and not from the integral quadrature.
Thus, residues and poles in the expansion~\eqref{eq:La_K_tilde} can be computed using various rational approximation algorithms, e.g., Pad\'e approximant~\cite{baker1996pade}, Best Uniform Rational Approximation~\cite{stahl2003best}, barycentric rational interpolation~\cite{berrut2005recent}, etc..
For more information on the rational interpolation methods, we refer to~\cite{trefethen2019approximation,celis2008practical}.
In this work, we employ the adaptive Antoulas--Anderson (AAA) algorithm~\cite[Fig. 4.1]{nakatsukasa2018aaa}, which in our test cases has demonstrated particular efficiency and robustness.



Let us consider the target function $f(z) = z^\alpha$ on the interval~$[a,b]$.
In our numerical simulation, we simply set $a$ to be equal to the time step size~$\dt$ of the time-integration scheme, and $b$ to the end time~$T$.
However, in the general case,  more sophisticated choices of the interval~$[a,b]$ can be considered, e.g., depending on the order of the scheme.
Following the AAA algorithm, the rational function is represented in barycentric form with interpolation at certain support points selected by the algorithm from a set of candidates provided by the user.
Since the target function grows faster near the origin,  we use a logarithmic grid on~$[\dt,T]$ as the candidate set.
Thus, we approximate the target function with a ratio of two polynomials of degree~$m$:
\begin{equation}\label{key}
	P(z)=\sum_{k=0}^{m}p_k z^k \qquad \text{and} \qquad Q(z)=\sum_{k=0}^{m}q_k z^k.
\end{equation}
Then, the rational function~$r(z) = P(z^{-1}) / Q(z^{-1}) = \sum_{k=0}^{m}p_k z^{m-k} / \sum_{k=0}^{m}q_k z^{m-k}$ approximates $f(z^{-1}) = z^{-\alpha}$ on the interval~$[\frac{1}{T},\frac{1}{h}]$.
Hence, the partial fractions decomposition of~$\hat{\tilde{K}}_{\alpha}(z) = r(z)$ yields the required multi-pole form~\eqref{eq:La_K_tilde}.
We apply the AAA algorithm to the reciprocal of the kernel spectrum, since it shows better stability for the values of $\alpha$ close to~$1$.

\begin{figure}[ht!]
	\centering\noindent
	\subfloat[]{
		\includegraphics[width=0.48\textwidth]{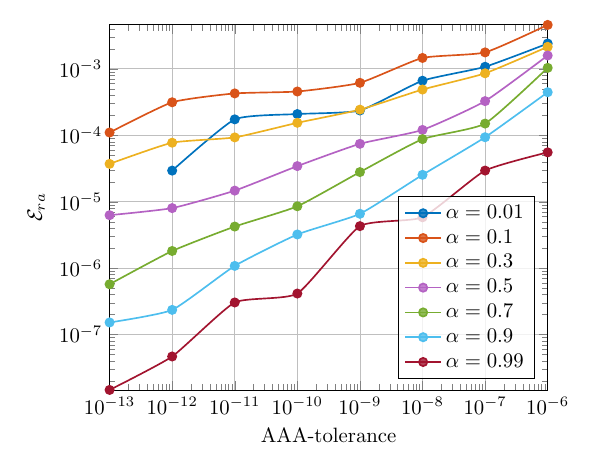}
		\label{fig:AAA_accuracy}
	}
	\hfill
	\subfloat[]{		
		\includegraphics[width=0.48\textwidth]{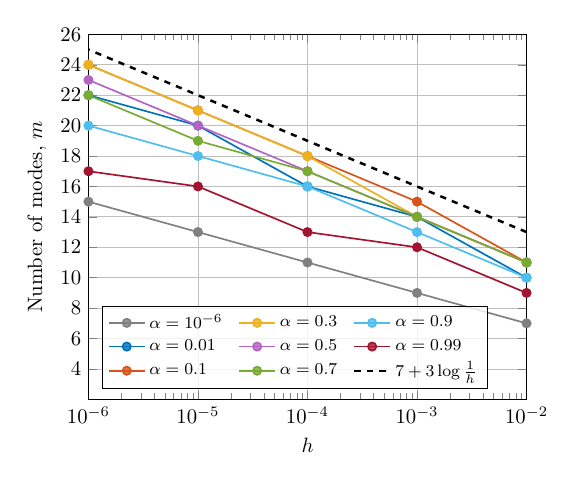}
		\label{fig:ModesNumber}
	}
	
	\caption{
		(a)~Dependence on the AAA algorithm tolerance of the rational approximation error~\eqref{eq:RA:error_def} for $T=1$ and $\dt=10^{-5}$.
		(b)~Dependence on~$h$ of the number of modes~$m$ with the AAA-tolerance $10^{-12}$, $T=1$.
		We observe that the number of modes $m(h)$ is below the heuristic bound $7 + 3\log \frac{1}{\dt}$.
	}
	\label{fig:RA_accuracy_Err}
\end{figure}

Since the Laplace transform is a compact operator from $\L{2}([0,\infty])$ to $\L{2}([0,\infty])$, the problem of its inversion on the real line is, generally speaking, ill-posed~\cite{epstein2008bad}.
Nevertheless, when the tolerance of the AAA algorithm is small enough, we observe convergence of the kernel approximation error, see, e.g.,~\Cref{fig:AAA_accuracy}.
Besides, approximation on the real line provides real polynomial coefficients $p_k$ and $q_k$.
To illustrate the accuracy, let us define the error of the kernel approximation as follows:
\begin{equation}\label{eq:RA:error_def}
\mathcal{E}_{ra} := \norm{K_{\alpha}-\tilde{K}^{Exp}_{\alpha}}_{\L{1}([h,T])} + \abs{\int_{0}^{h}\left(K_{\alpha}(s)-\tilde{K}_{\alpha}^{Exp}(s)\right)\d s-w_{\infty}},
\end{equation}
where we denote the sum-of-exponentials part of the kernel approximant by
\begin{equation}
\tilde{K}_{\alpha}^{Exp}(t) := \sum_{k=1}^{m} w_k\exp{-\lambda_k t}.
\end{equation}
Thus, we split the error into two terms: the error of the sum-of-exponential approximation and the error of the local contribution due to the $w_{\infty}$ term (see also the error estimate in~\Cref{th:error_estimate}).
In~\Cref{fig:AAA_accuracy}, we show the convergence of this error with respect to the tolerance of the AAA algorithm for $T=1$ and $\dt=10^{-5}$.
The integral in~\eqref{eq:RA:error_def} is computed using $scipy$ package~\cite{scipy}.

In~\Cref{fig:ModesNumber}, the number of modes~$m$ with the AAA-tolerance $10^{-12}$ is shown as function of the step size~$h$ for different values of the fractional power~$\alpha$.
Without loss of generality, we fix the time interval size $T=1$, since increasing the time interval with the same number of steps can be seen as decreasing~$h$ for the rescaled time ${\tilde{t}=t/T}$.
We observe that the number of modes $m$ grows as $\log N$, which is typical for the majority of kernel compression methods.
However, we remark that the value itself of the number of modes is significantly smaller (order of $20$) in comparison with many other methods, where the number of modes (auxiliary variables, quadrature points) is typically of order of hundreds, see, e.g.,~\cite{baffet2017kernel,baffet2019gauss,jiang2017fast,li2010fast,zeng2018stable}.
We can also observe that the number of modes decreases at the limits of the interval $(0,1)$.
Besides, for the value~$\alpha=1$, we obtain $m=1$ with the single mode $w_1=1$ and $\lambda_1=0$.
For $\alpha=0$, we have $m=0$, when only the "infinity" mode $w_{\infty}=1$ remains.
Note that the coefficient~$w_{\infty}$ takes values between $0$ and $1$, where the extremities correspond respectively to $\alpha=1$ and $0$.

\vspace{-3ex}

\begin{figure}[ht!]
	\centering\noindent
	
	\subfloat[]{
		\includegraphics[width=0.45\textwidth]{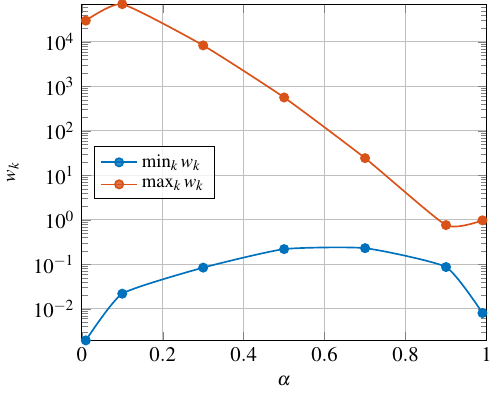}
		\label{fig:weights_bounds}
	}
	\hfill
	\subfloat[]{		
		\includegraphics[width=0.45\textwidth]{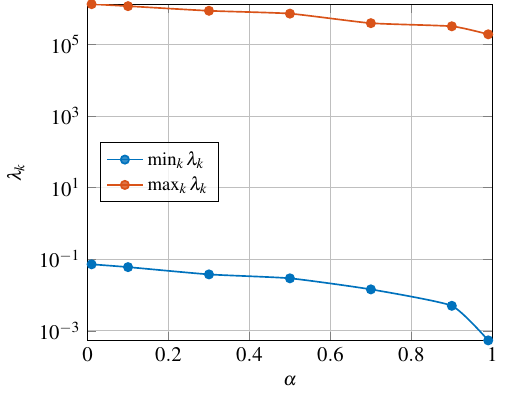}
		\label{fig:poles_bounds}
	}
	
	\caption{
		The maximum and the minimum values of the weights $w_k$ (left) and the exponents $\lambda_k$ (right) obtained using AAA rational approximation as functions of the factional power $\alpha$.
		The AAA tolerance is fixed to $10^{-12}$, $h=10^{-5}$ and $T=1$.
	}
	\label{fig:weights_poles_bounds}
\end{figure}

\begin{figure}[ht!]
	\centering\noindent
	
	\includegraphics[width=0.5\textwidth]{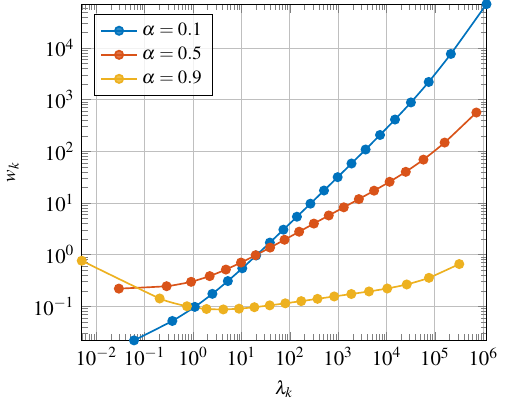}
	
	\caption{
		The distribution of the weights $w_k$ and the poles $\lambda_k$ for $\alpha=0.1, 0.5, 0.9$.
		The nodes correspond to the actual values, while the line connectors are for illustration purpose.
	}
	\label{fig:weights_poles_distribution}
\end{figure}

The kernel $\FrKer{\alpha}(t)=t^{\alpha-1}/\Gamma(\alpha)$, $\alpha\in(0,1)$, is a completely monotone function.
Besides, its sum-of-exponentials approximation $\tilde{K}_{\alpha}^{Exp}(t) := \sum_{k=1}^{m} w_k\exp{-\lambda_k t}$ is also completely monotone if the weights~${w_k}$ and the exponents~${\lambda_k}$ are positive, according to Berstein theorem on completely monotone functions~\cite{bernstein1929fonctions,widder2015laplace}; see also~\cite{kammler1977prony,braess2012nonlinear}.
Note that for the results presented in~\Cref{fig:AAA_accuracy}, the AAA algorithm provides $w_k\ge 0$ and~$\lambda_k\ge 0$, and therefore the resulting sum-of-exponentials approximation is completely monotone. 
The maximum and the minimum values of $w_k$ and $\lambda_k$ are depicted in~\Cref{fig:weights_poles_bounds} for different values of $\alpha$ (same as in~\Cref{fig:AAA_accuracy}), the AAA tolerance $10^{-12}$, $h=10^{-5}$ and $T=1$.
It shows that the weights and the exponents stay positive.
Thus, the AAA algorithm shows excellent robustness with respect to a large range of $\alpha$ values and a small tolerance. 
We note that negative spurious poles only occur for very small $\alpha$ with very small AAA tolerance. However, in many applications the relevant $\alpha$ does not approach zero; see, e.g.,~\cite{schmidt2002finite,meral2010fractional} in viscoelasticity, \cite{valentim2020can} in tumor growth modeling. 
In addition, we illustrate in~\Cref{fig:weights_poles_distribution} the distribution of the weights $w_k$ and the poles $\lambda_k$ for $\alpha=0.1, 0.5, 0.9$.

	\section{Error analysis}
	\label{sec:EA}

In this section, we estimate the bound for the global error on the interval~$[0,T]$ between the solution to~\eqref{eq:def:solution} and its approximation~\eqref{eq:RA_expansion}.
The bound is suggested in~\Cref{th:error_estimate} below.
Beforehand, we introduce the following auxiliary lemmas required for the proof of the theorem.


\begin{lemma}[Gr\"onwall inequality~\cite{ye2007generalized}]\label{th:Gronwall}
	Let $v, \eps\in\L{2}([0,T])$ and a constant $C>0$ be such that for all $t\in[0,T]$
	\begin{equation}\label{eq:Gronwall:given}
	\abs{v(t)}
	\le	
	\abs{\eps(t)}
	+
	C\,\left(
	K_{\alpha}\conv\abs{v}
	\right)(t).
	\end{equation}
	Then, $v(t)$ satisfies the following bound:
	\begin{equation}\label{key}
	\abs{v(t)}
	\le	
	\abs{\eps(t)}
	+
	C\int_{0}^{t}s^{\alpha- 1}\Ea_{\alpha,\alpha}[C s^{\alpha}]\cdot\abs{\eps(t-s)} \d s,
	\qquad
	t\in[0,T],
	\end{equation}
	where $\Ea_{\alpha,\beta}[x]$ denotes the Mittag-Leffler function~\cite{Kexue2011}:
	\begin{equation}\label{eq:Mittag-Leffler}
	\Ea_{\alpha,\beta}[x] := \sum_{k=0}^{\infty}\frac{x^k}{\Gamma(\alpha k + \beta)}.
	\end{equation}
\end{lemma}


\begin{lemma}\label{th:aux1}	
	For $\alpha\in(0,1]$ and $C>0$, it holds that
	\begin{equation}\label{key}
	K_{\alpha}(t) 	+
	C\int_{0}^{t}s^{\alpha- 1}\Ea_{\alpha,\alpha}[C s^{\alpha}]\,K_{\alpha}(t-s)  \d s
	=
	t^{\alpha-1} \Ea_{\alpha,\alpha}[Ct^{\alpha}].
	\end{equation}	
\end{lemma}

\begin{proof}
	Note that for any $\nu\ge 0$, the fractional integral of the function $t^{\nu-1}$ can be computed by the following formula, see~\cite{kilbas2006theory,pang2018fractional},
	\begin{equation}
	\FrInt{\alpha} t^{\nu-1} = \frac{\Gamma(\nu)}{\Gamma(\nu+\alpha)} t^{\nu+\alpha-1}.
	\end{equation}
	Using this, we can then write
	\begin{align}\label{key}
	C\int_{0}^{t}s^{\alpha- 1}\Ea_{\alpha,\alpha}[C s^{\alpha}]\,K_{\alpha}(t-s)  \d s
	&=
	\FrInt{\alpha}\sum_{k=1}^{\infty}\frac{C^k\,t^{\alpha k - 1}}{\Gamma(\alpha k)}
	=
	\sum_{k=1}^{\infty}\frac{C^k\,t^{\alpha k + \alpha - 1}}{\Gamma(\alpha k+\alpha)}
	\\
	&=
	t^{\alpha-1} \left(\Ea_{\alpha,\alpha}[Ct^{\alpha}] - \frac{1}{\Gamma(\alpha)}\right).
	\end{align}
	Hence, it follows the statement of the lemma.
\end{proof}


\vspace{3ex}
It is known that the solution of a fractional differential equations often exhibits a weak singularity at the initial time~\cite{lubich1986discretized}.
The following lemma estimates the asymptotic behavior of the derivative.

\begin{lemma}\label{th:difference}
	Let $u(t) = u_0 + K_{\alpha}\conv F[u]$ be defined as in~\Cref{sec:Preliminaries}.
	And let $F: \Bochner \to \BochnerDual$ be a uniformly Lipschitz continuous operator, i.e., there exist $C_1>0$ such that for any $v_1,v_2\in\Bochner$, $v_3\in\HS$ and $t\in[0,T]$,
	\begin{equation}\label{eq:EA:continuityF}
	\abs{\dualH{F[v_1](t) - F[v_2](t)}{v_3}} \le C_1 \normH{v_1(t) - v_2(t)}\cdot\normH{v_3}.
	\end{equation}
	Then, for $t>0$ and $h>0$, the difference $u(t+\dt)-u(t)$ is bounded by
	\begin{equation}\label{key}
	\normH{u(t+\dt)-u(t)}
	\le {\dt}\,t^{\alpha-1} \sup\limits_{[0,\dt]}\normHdual{F[u]}\,
	 \Ea_{\alpha,\alpha}[Ct^{\alpha}].
	\end{equation}
\end{lemma}

\begin{proof}
	By definition of $u$, we can write
	\begin{align}\label{key}
	u(t+\dt)-u(t)
	=&
	K_{\alpha}\conv F[u](t+\dt) - K_{\alpha}\conv F[u](t) \\
	=&
	\int_{0}^{t} K_{\alpha}(t-s)\,\left[F[u](s+\dt) - F[u](s)\right]\d s\\
	&+ \int_{0}^{\dt} K_{\alpha}(t+\dt-s)\,F[u](s) \d s.
	\end{align}
	Let us denote the difference by $v(t):=\normH{u(t+\dt)-u(t)}$.
	Then, from the previous equation, we have
	\begin{align}\label{key}
	v(t)
	\le&
	C_1\int_{0}^{t} K_{\alpha}(t-s)\,v(s)\d s
	+ K_{\alpha}(t)\int_{0}^{\dt} \left(1+\frac{\dt-s}{t}\right)^{\alpha-1}\,\normHdual{F[u](s)} \d s \\	
	\le&
	C_1\int_{0}^{t} K_{\alpha}(t-s)\,v(s)\d s
	+ {\dt}\,K_{\alpha}(t) \sup\limits_{[0,\dt]}\normHdual{F[u]}.
	\end{align}
	Hence, by~\Cref{th:Gronwall} and~\Cref{th:aux1}, we obtain
	\begin{align}\label{key}
	v(t)
	&\le
	{\dt}\, \sup\limits_{[0,\dt]}\normHdual{F[u]}
	\left(K_{\alpha}(t)
	+
	C_1\int_{0}^{t}s^{\alpha- 1}\Ea_{\alpha,\alpha}[C_1 s^{\alpha}]\,K_{\alpha}(t-s) \d s
	\right) \\
	&\le
	{\dt}\, \sup\limits_{[0,\dt]}\normHdual{F[u]}\,
	t^{\alpha-1} \Ea_{\alpha,\alpha}[C_1 t^{\alpha}].
	\end{align}
	
\end{proof}


	\vspace{3ex}
	Thus, the derivative of the solution to~\eqref{eq:def:solution} is bounded by $\O(t^{\alpha-1})$.
	The weak singularity at the initial time can pollute the local error at the beginning of the time interval~\cite{diethelm2020good}.
	So, the optimal convergence rate can be not observed globally using classical Lebesgue norms in time. 
	Therefore, we introduce the weighted Bochner space~$\Bochner_{\alpha}$ defined by the following norm:
\begin{equation}\label{key}
	\normBA{u} = \normB{u\cdot t^{1-\alpha}} = \left(\int_{0}^{T}\normH{u(t)}^2\, t^{2\cdot(1-\alpha)} \d t\right)^{\frac{1}{2}}.
\end{equation}
We are now ready to prove the following theorem providing an error bound for the approximate solution~$\tilde{u}$.

\begin{theorem}\label{th:error_estimate}
	Let $F: \Bochner \to \BochnerDual$ be a uniformly Lipschitz continuous operator, i.e., there exist $C_1>0$ such that \eqref{eq:EA:continuityF} holds.
	And let
	\begin{equation}
	u(t) = u_0 + K_{\alpha}\conv F[u],
	\qquad
	\tilde{u}(t) =  u_0 + \tilde{K}_{\alpha}\conv F[\tilde{u}],
	\end{equation}
	such that for $h>0$, there exists a constant $C_2>0$, depending on $\alpha$, that
	\begin{equation}\label{eq:EA:statement_bounds}
	\abs{\int_{0}^{h}\left(K_{\alpha}(s)-\tilde{K}_{\alpha}^{Exp}(s)\right)\d s-w_{\infty}} + \norm{K_{\alpha}-\tilde{K}_{\alpha}^{Exp}}_{\L{1}([h,T])} \le C_2\,h^{1+\alpha},
	\end{equation}
	and $\abs{K_{\alpha}(t)-\tilde{K}_{\alpha}^{Exp}(t)}\le K_{\alpha}(t)$ for $t\in[0,\dt]$.
	Then, the following error estimate holds:
	\begin{equation}\label{eq:ErrorBound}
	\norm{u - \tilde{u}}_{\Bochner_{\alpha}}
	\le 
	C(\alpha,T)\, \dt^{1+\alpha}\,\left(\sup\limits_{[0,\dt]}\normHdual{\drv{t}{^\alpha u}} + \normBHdual{\drv{t}{^\alpha u}}\right),
	\end{equation}
	where the constant $C(\alpha,T)$ depends only on $\alpha$ and $T$.
\end{theorem}

\begin{proof}
	Let us introduce $\tilde{v}(t) := u_0 + \tilde{K}_{\alpha}\conv F[u](t)$.
	Then, adding $0=\tilde{v}-\tilde{v}$, we find for the norm
	\begin{equation}\label{eq:EA:compl_term}
	\normH{u - \tilde{u}}^2
	\le
	\abs{\scaH{u - \tilde{v}}{u - \tilde{u}}}
	+
	\abs{\scaH{\tilde{v} - \tilde{u}}{u - \tilde{u}}}.	
	\end{equation}
	We denote by $g := K_{\alpha}-\tilde{K}_{\alpha}^{Exp}$ the difference of the kernels, and by $g_{[a,b]}$ a function which coincides with $g$ on $[a,b]$ and vanishes elsewhere.
	Then, we can formally write
	\begin{align}\label{key}
	(K_{\alpha}-\tilde{K}_{\alpha})\conv F[u](t) 
	&=
	\int\limits_{0}^{\mathclap{\min\{h,t\}}} g(s) F[u](t-s) \d s - w_{\infty}\,F[u](t) 
	+ g_{[h,T]} \conv F[u](t) .
	\end{align}
	Adding $0=F[u](t)\int_{0}^{h} g(s) \d s - F[u](t)\int_{0}^{h} g(s) \d s$, this can be rewritten as
	\begin{align}\label{key}
	(K_{\alpha}-\tilde{K}_{\alpha})\conv F[u](t)
	=&
	\int\limits_{0}^{\mathclap{\min\{h,t\}}} g(s) \left[F[u](t-s)-F[u](t)\right] \d s - \int\limits_{\mathclap{\min\{h,t\}}}^{h} g(s) \d s \cdot F[u](t) + \\
	&+ \left[\int_{0}^{h} g(s) \d s - w_{\infty}\right] F[u](t) 
	+ g_{[h,T]} \conv F[u](t).
	\end{align}
	Then, using the continuity assumption~\eqref{eq:EA:continuityF}, we can bound the first term in~\eqref{eq:EA:compl_term} as follows:
	\begin{align}\label{eq:EA:1nd_term}
	\abs{\scaH{u - \tilde{v}}{u - \tilde{u}}}
	&= \abs{\scaH{(K_{\alpha}-\tilde{K}_{\alpha})\conv F[u]}{u - \tilde{u}}} \\
	&\le
	\left(\eps_1(t) + \eps_2(t)\right)\cdot\normH{u - \tilde{u}},
	\end{align}
	where the terms $\eps_1(t)$ and $\eps_2(t)$ are respectively defined as
	\begin{align}
		\eps_1(t)
		&:=
		C_1\int\limits_{0}^{\mathclap{\min\{h,t\}}} \abs{g(s)}\cdot \normH{u(t-s)-u(t)} \d s +
		 \frac{\dt^{\alpha}\,\Ind_{[0,h]}(t)}{\Gamma(\alpha+1)}\,\normHdual{F[u]},
		\label{eq:eps_def1} \\
		\eps_2(t)
		&:= \abs{\int_{0}^{\dt} g(s) \d s - w_{\infty}}\cdot\normHdual{F[u]} + \abs{g_{[\dt,T]}} \conv \normHdual{F[u]}.
		\label{eq:eps_def2}		
	\end{align}
	Let us investigate the term~$\eps_1(t)$.
	First, using \Cref{th:difference}, we can write 
	\begin{equation}
	\int\limits_{0}^{\mathclap{\min\{h,t\}}} \abs{g(s)}\cdot
	\normH{u(t-s)-u(t)} \d s
	\le 
	\sup\limits_{[0,\dt]}\normHdual{F[u]}\Ea_{\alpha,\alpha}[C_1 T^{\alpha}]\,\int\limits_{0}^{\mathclap{\min\{h,t\}}} s K_{\alpha}(s)\, 
	(t-s)^{\alpha-1} \d s.
	\label{eq:est1}
	\end{equation}
	Note that given $t \ge t_0>0$ and $\alpha\in(0,1]$, it holds that
	\begin{equation}\label{eq:aux2}
	\int_{0}^{t_0} (t-s)^{\alpha-1}\d s
	=
	t^{\alpha-1}\, t_0\int_{0}^{1} \left(1-s\frac{t_0}{t}\right)^{\alpha-1}\d s
	\le 
	t^{\alpha-1}\, t_0\int_{0}^{1} \left(1-s\right)^{\alpha-1}\d s
	=
	t^{\alpha-1}\, t_0 / \alpha.
	\end{equation}	
	Then, the integral in the right hand side of~\eqref{eq:est1} can be bounded using~\eqref{eq:aux2} with $t_0=\min\{h,t\}$:
	\begin{equation}\label{eq:est2}
	\int\limits_{0}^{\mathclap{\min\{h,t\}}} s K_{\alpha}(s)\, 
	(t-s)^{\alpha-1} \d s
	\le
	\frac{\dt^{\alpha}}{\Gamma(\alpha)}\,\int\limits_{0}^{\mathclap{\min\{h,t\}}}(t-s)^{\alpha-1} \d s
	\le	
	\frac{\dt^{1+\alpha} t^{\alpha-1}}{\Gamma(\alpha+1)}.
	\end{equation}
	Substituting \eqref{eq:est1} and \eqref{eq:est2} to \eqref{eq:eps_def1}, we obtain the bound for $\eps_1(t)$:
	\begin{equation}\label{eq:eps1_bound}
		\eps_1(t) \le C_3 \,\dt^{1+\alpha}K_{\alpha}(t)\,\sup\limits_{[0,\dt]}\normHdual{F[u]}
		 + \frac{\dt^{\alpha}\,\Ind_{[0,h]}(t)}{\Gamma(\alpha+1)}\,\normHdual{F[u]},
	\end{equation}
	where we denoted $C_3:=\Ea_{\alpha,\alpha}[C_1 T^{\alpha}] \, \Gamma(\alpha)/\Gamma(\alpha+1)$.
	
	For the second term in~\eqref{eq:EA:compl_term}, we also use the continuity~\eqref{eq:EA:continuityF} to obtain the following upper bound:
	\begin{align}\label{eq:EA:2nd_term}
	\abs{\scaH{\tilde{v} - \tilde{u}}{u - \tilde{u}}}
	&= \abs{\scaH{K_{\alpha}\conv (F[u] - F[\tilde{u}])}{u - \tilde{u}}}
	\\
	&\le
	C_1\left(K_{\alpha}\conv\normH{u - \tilde{u}}\right)\cdot\normH{u - \tilde{u}}.
	\end{align}
	Then, substituting~\eqref{eq:EA:1nd_term} and~\eqref{eq:EA:2nd_term} to~\eqref{eq:EA:compl_term}, we obtain
	\begin{equation}
	\normH{u - \tilde{u}}
	\le	
	\eps_1(t) + \eps_2(t)
	+
	C_1\left(
	K_{\alpha}\conv\normH{u - \tilde{u}}
	\right).
	\end{equation}
	Hence, by the Gr\"onwall inequality (\Cref{th:Gronwall}), we have
	\begin{equation}\label{eq:final_bound1}
	\normH{u - \tilde{u}}
	\le	
	\eps_1(t) + \eps_2(t)
	+
	C_1\int_{0}^{t}(t-s)^{\alpha- 1}\Ea_{\alpha,\alpha}[C_1 (t-s)^{\alpha}]\,\left(\eps_1(s) + \eps_2(s)\right) \d s.
	\end{equation}	
	Moreover, using the bound~\eqref{eq:eps1_bound} and applying~\Cref{th:aux1} and~\eqref{eq:aux2}, we obtain  
	the following inequality:
	\begin{equation}\label{key}	
	\begin{split}
	\eps_1(t)
	+
	C_1&\int_{0}^{t}(t-s)^{\alpha- 1}\Ea_{\alpha,\alpha}[C_1 (t-s)^{\alpha}]\,\eps_1(s)\d s
	\le 
	\frac{\dt^{\alpha}\,\Ind_{[0,h]}(t)}{\Gamma(\alpha+1)}\normHdual{F[u]}
	\\
	&
	+ \left(
	C_3\,
	\dt^{1+\alpha}\, t^{\alpha-1}\Ea_{\alpha,\alpha}[C_1 T^{\alpha}]
	+
	C_1\,\frac{\dt^{1+\alpha}t^{\alpha-1}}{\alpha \Gamma(1+\alpha)}\Ea_{\alpha,\alpha}[C_1 T^{\alpha}]
	\right)\,
	\sup\limits_{[0,\dt]}\normHdual{F[u]}.
	\end{split}
	\end{equation}
	Thus, \eqref{eq:final_bound1} reads
	\begin{equation}\label{key}
	\begin{split}
	\normH{u - \tilde{u}}
	\le	 C_4\dt^{1+\alpha}\,t^{\alpha-1}\,\sup\limits_{[0,\dt]}\normHdual{F[u]}
	+ \frac{\dt^{\alpha}\,\Ind_{[0,h]}(t)}{\Gamma(\alpha+1)}\normHdual{F[u]}
	\\
	+ \eps_2(t) + C_1\int_{0}^{t}(t-s)^{\alpha- 1}\Ea_{\alpha,\alpha}[C (t-s)^{\alpha}]\,\eps_2(s) \d s,
	\end{split}
	\end{equation}
	where we denoted $C_4:=C_3\left(\Ea_{\alpha,\alpha}[C_1 T^{\alpha}] + C_1/\Gamma(\alpha+1)\right)$.
	Then, computing the $\L{2}$-norm, weighted with $t^{2\cdot(1-\alpha)}$, we obtain that
	\begin{equation}\label{eq:error_estimate_1}
	\norm{u - \tilde{u}}_{\Bochner_{\alpha}}
	\le 
	C_4 \dt^{1+\alpha}\,\sqrt{T}\,\sup\limits_{[0,\dt]}\normHdual{F[u]}
	+
	\frac{\dt^{1+\alpha}\,\normBHdual{F[u]}}{\Gamma(\alpha+1)}
	+
	T^{1-\alpha} \Ea_{\alpha,1}[C_1\,T^\alpha]\cdot\norm{\eps_2}_{\L{2}([0,T])}
	,
	\end{equation}
	where we used Young's convolution inequality for the last two terms.
	Applying again Young's convolution inequality, along with the assumption~\eqref{eq:EA:statement_bounds}, we thus obtain from~\eqref{eq:eps_def2}:
	\begin{equation}
	\norm{\eps_2}_{\L{2}([0,T])} 
	\le
	\left(\abs{\int_{0}^{h} g(s) \d s - w_{\infty}} + \norm{g}_{\L{1}([h,T])}\right)\cdot\normBHdual{F[u]}
	\le C_2\,h^{1+\alpha}\,\normBHdual{\drv{t}{^\alpha u}}.	\label{eq:norm_eps}
	\end{equation}
	Then, substituting~\eqref{eq:norm_eps} to \eqref{eq:error_estimate_1},
	we eventually obtain the error estimate
	\begin{equation}\label{key}
	\norm{u - \tilde{u}}_{\Bochner_{\alpha}}
	\le 
	C(\alpha,T)\; \dt^{1+\alpha}\,\left(\sup\limits_{[0,\dt]}\normHdual{\drv{t}{^\alpha u}} + \normBHdual{\drv{t}{^\alpha u}}\right),
	\end{equation}
	where the constant $C$ depends on $\alpha$ and $T$.
	$ $
\end{proof}

\begin{figure}[ht!]
	\centering\noindent
	
	\includegraphics[width=0.5\textwidth]{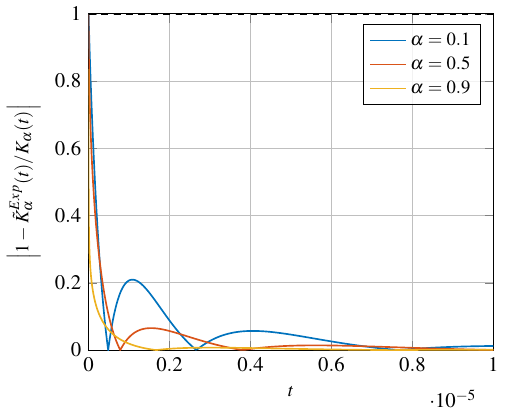}
	
	\caption{
		Numerical illustration that the assumption $|K_{\alpha}(t)-\tilde{K}_{\alpha}^{Exp}(t)|\le K_{\alpha}(t)$, $t\in[0,\dt]$, in~\Cref{th:error_estimate} holds. Here, $\alpha=0.1, 0.5, 0.9$, the AAA tolerance $10^{-12}$, $h=10^{-5}$ and $T=1$.
	}
	\label{fig:misfit_near_origin}
\end{figure}

Note that the bound~\eqref{eq:EA:statement_bounds} is illustrated in~\Cref{fig:AAA_accuracy}.
The assumption $\abs{K_{\alpha}(t)-\tilde{K}_{\alpha}^{Exp}(t)}\le K_{\alpha}(t)$ for $t\in[0,\dt]$ is easy to satisfy, since $K_{\alpha}(t)$ explodes at the origin.
For an illustration, we also show in~\Cref{fig:misfit_near_origin} the function $|1-\tilde{K}_{\alpha}^{Exp}(t)/K_{\alpha}(t)|$ on the interval $(0,\dt]$ for $\alpha=0.1, 0.5, 0.9$, the AAA tolerance $10^{-12}$, $h=10^{-5}$ and $T=1$.

	\section{Numerical schemes}
	\label{sec:Scheme}

\label{sec:Schemes}


\newcommand{\uhk}{u_{h,k}}
\newcommand{\uht}{\tilde{u}_h}
\newcommand{\errh}{\epsilon_h}
\newcommand{\errhk}{\epsilon_{h,k}}


Let us denote by~$\uht$ the numerical approximation of~$\tilde{u}$, defined by~\eqref{eq:RA_expansion}, and by~$\uhk$ the approximations of the modes~$u_k$, $k=1,\ldots,m$, in~\eqref{eq:sys_modes}.
In what follows, we use the superscript $n$ to indicate the time step~$n$.
In particular, we define $t^n = n\,\dt$, where $h$ stands for the time step size, and $\uht^{n}=\uht(t^{n})$.
For the sake of simplicity, we will use the notation $F[\uht^{n}]$ to denote $F[\uht](t^n)$.
The discretized solution of the system of equations~\eqref{eq:sys_modes} can be numerically computed with any suitable numerical scheme.
The simplest case of the so-called $\theta$-scheme, including Euler and Crank-Nicolson time-integration, is introduced in the following proposition.


\begin{proposition}\label{th:scheme}
	Applying a standard $\theta$-scheme to the modal system~\eqref{eq:sys_modes}, it yields the following scheme for~\eqref{eq:RA_expansion}: 
	\begin{equation}\label{eq:final_scheme}
	\uht^{n+1} - \beta^1\,F[\uht^{n+1}] = \beta^2\, F[\uht^{n}] + u_0 + \sum_{k=1}^{m}\gamma_k\,\uhk^{n},
	\end{equation}
	where the discrete modes $\uhk^{n}$, $k=1,\ldots,m$, are updated by
	\begin{equation}\label{eq:final_scheme_modes}
	\begin{aligned}	
	\uhk^{n+1} = \gamma_k\,\uhk^{n} + \beta_k^1\,F[\uht^{n+1}] + \beta_k^2\,F[\uht^{n}],
	\qquad
	\uhk^{0} = 0,
	\end{aligned}
	\end{equation}
	with coefficients
	\begin{equation}\label{eq:scheme:coefs}	
	\begin{gathered}
	\beta^1_k = \frac{w_k\theta\dt}{1 + \theta \lambda_k\dt}, \qquad
	\beta^2_k = \frac{w_k(1-\theta)\dt}{1 + \theta \lambda_k\dt}, \qquad
	\gamma_k = \frac{1-(1-\theta)\lambda_k \dt}{1 + \theta \lambda_k\dt},
	\\	
	\beta^1 = \sum_{k=1}^{m}\beta^1_k + w_{\infty}, 
	\qquad 
	\beta^2 = \sum_{k=1}^{m}\beta^2_k.
	\end{gathered}
	\end{equation}
\end{proposition}

\begin{remark}
	Let us note that the coefficient~$\beta^1$ is nothing else than the rational approximation~\eqref{eq:La_K_tilde} of $(\theta\dt)^\alpha$.
	In particular, we have $\beta^1=(\theta\dt)^\alpha$ if $z=\theta\dt$ is a support point of the rational approximation.
\end{remark}

\begin{remark}
	When $\alpha=1$, there is only one mode $w_1=1$, $\lambda_1=0$, therefore, the above scheme reduces to the classical integer-order $\theta$-scheme.
\end{remark}

\begin{figure}[ht!]
	\centering\noindent
	\includegraphics[width=\textwidth]{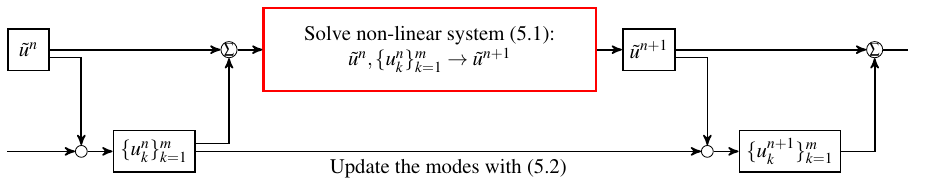}
	\caption{
		Scheme of the $n$-th time iteration, representing the system~\eqref{eq:final_scheme}-\eqref{eq:final_scheme_modes}.
		Both the integral approximation~$\tilde{u}^n$ and the family of modes~$\{u_k^{n}\}_{k=0}^{m}$ have to be updated in each time step.
		However, each step requires only one non-linear system of the original size to solve.
	}
	\label{fig:TimeStepScheme}
\end{figure}

Let us note that the proposed scheme does not require the solution of a large coupled system of equations, but consists in alternating updates of the full solution~$\uht$, Eq.~\eqref{eq:final_scheme}, and updates of the modes~$\uhk$, Eq.~\eqref{eq:final_scheme_modes}.
A graphical illustration of the algorithm is suggested in~\Cref{fig:TimeStepScheme}.	
The modes updates~\eqref{eq:scheme:coefs} are completely decoupled and can be computed in parallel.
Besides, they are linear.
A non-linear equation of the original size has to be solved only once per time-step in~\eqref{eq:final_scheme}, using any preferred non-linear solver (e.g., Newton-Raphson or Fixed-point).
In particular, in the PDE case, when $F$ involves a spacial differential operator, the PDE system is solved only in~\eqref{eq:final_scheme}.
Moreover, for the updates~\eqref{eq:scheme:coefs}, one does not even have to solve a mass matrix system.
Indeed, instead of computing the modes $\uhk^{n}$, $k=1,\ldots,m$, themselves, one can proceed with numerical integration computing only $M\uhk^{n}$, where $M$ stands for the formal mass matrix.
And no explicit computation of the modes is necessary for computing the full solution~$\uht^{n}$.
	
The proposed $\theta$-scheme is simple, however, it does not guarantee unconditional stability for an arbitrary operator~$F$.
In particular, unconditionally stable schemes are usually based on the splitting of the operator~\cite{eyre1998unconditionally,eyre1998unconditionally2}.
So, let us consider $F$ in the form $F[v] = F_{-}[v] + F_{+}[v]$, where the monotonous operators $F_{-}$ and $F_{+}$ corresponds to the decreasing and strictly increasing parts of $F$, respectively.
That is, for all $v_1,v_2\in\Bochner$, it holds
\begin{equation}\label{key}
	\scaH{F_{-}[v_1]-F_{-}[v_2]}{v_1-v_2} \le 0,
	\qquad
	\scaH{F_{+}[v_1]-F_{+}[v_2]}{v_1-v_2} > 0.
\end{equation}
In addition, we rewrite the continuity condition with some~$C_F>0$ in the form
\begin{equation}\label{key}
	\normHdual{F_{-}[v_1]-F_{-}[v_2]} + \normHdual{F_{+}[v_1]-F_{+}[v_2]} \le C_F\normH{v_1-v_2}.
\end{equation}
In the following lemmas, we propose numerical schemes for such $F$ and estimate the associated discretization error $\errh^{n}:=\normH{\tilde{u}(t^{n}) - \uht^{n}}$.
We also introduce the modal discretization errors $\errhk^{n}:=\normH{u_k(t^{n}) - \uhk^{n}}$.



\begin{lemma}[Implicit Euler]\label{th:IE}
	Let $\uht^{n}$ be defined by the following time-stepping scheme:
	\begin{equation}\label{eq:IE:scheme}
	\uht^{n+1} - \beta\,F_{-}[\uht^{n+1}] = \beta\,F_{+}[\uht^{n}] + u_0 + \sum_{k=1}^{m}\gamma_k\,\uhk^{n},
	\end{equation}
	where the discrete modes $\uhk^{n}$, $k=1,\ldots,m$, are updated by
	\begin{equation}\label{eq:IE:modes}
	\begin{aligned}	
	\uhk^{n+1} = \gamma_k\,\uhk^{n} + \beta_k\,F[\uht^{n+1}],
	\qquad
	\uhk^{0} = 0,
	\end{aligned}
	\end{equation}
	with coefficients given by
	\begin{equation}\label{eq:IE:coefs}	
	\gamma_k = \frac{1}{1 + \lambda_k\dt}, \qquad
	\beta_k = \frac{w_k\dt}{1 + \lambda_k\dt}, \qquad
	\beta = \sum_{k=1}^{m}\beta_k + w_{\infty}.
	\end{equation}
	Then, $\uht^{n}$ approximates $\tilde{u}(t^n)$ with the discretization error of order~$\dt$,
	\begin{equation}\label{key}
	\errh^{n} = \normH{\tilde{u}(t^{n}) - \uht^{n}} \le \O(m\dt) \cdot \exp{C_F \left[(t^n)^\alpha+\eps_{ra}\right]},
	\end{equation}
	where $\eps_{ra}:=\max\limits_{s\in[\frac{1}{T},\frac{1}{\dt}]}\left|\hat{K}(s)-\hat{\tilde{K}}(s)\right|$ stands for the rational approximation error.
\end{lemma}

\begin{proof}	
	Taylor expansion of $u_k(t^{n})$ at the point~$t^{n+1}$, with the first derivative given by~\eqref{eq:sys_modes}, yields
	\begin{equation}
		u_k(t^{n+1}) = u_k(t^{n}) - \dt\,\underbrace{\left(\lambda_k\,u_k(t^{n+1}) - w_k F[\tilde{u}](t^{n+1})\right)}_{\drv{t}{u_k}(t^{n+1})} + \O(\dt^2).
	\end{equation}
	Hence, we express $u_k(t^{n+1})$:
	\begin{equation}\label{eq:conv:IE:uhk}
		u_k(t^{n+1}) = \gamma_k u_k(t^{n}) + \beta_k\,F[\tilde{u}](t^{n+1}) + \gamma_k\O(\dt^2).
	\end{equation}
	Summing up the modes, $u_0$ and $u_{\infty}$, we obtain
	\begin{equation}
		\tilde{u}(t^{n+1}) - \beta\,F[\tilde{u}](t^{n+1}) = u_0 + \sum_{k=1}^{m}\gamma_k u_k(t^{n}) + \O(m\dt^2).
	\end{equation}
	Note that $\beta=\O(h^\alpha)$.
	Then, using Taylor expansion of $F_{+}[\tilde{u}](t^{n+1})$ at the point $t^{n}$, we can write
	\begin{equation}\label{eq:conv:IE:uh}
		\tilde{u}(t^{n+1}) - \beta\,F_{-}[\tilde{u}](t^{n+1}) = \beta\,F_{+}[\tilde{u}](t^{n}) + u_0 + \sum_{k=1}^{m}\gamma_k u_k(t^{n}) + \O(m\dt^2+h^{1+\alpha}).
	\end{equation}
	Recall that $F$ is Lipschitz continuous and $F_{-}$ monotonously decreases, which implies
	\begin{equation}\label{eq:proof:continuity_and_monotonicity}
	\begin{aligned}
	\abs{\scaH{F[\tilde{u}](t^{n})-F[\uht]}{\tilde{u}(t^{n}) - \uht^n}} &\le C_F \abs{\errh^{n}}^2,
	\\
	\scaH{F_{-}[\tilde{u}](t^{n})-F_{-}[\uht]}{\tilde{u}(t^{n}) - \uht^n} &\le 0.
	\end{aligned}		
	\end{equation}
	Thus, subtracting~\eqref{eq:IE:scheme} and~\eqref{eq:IE:modes} from~\eqref{eq:conv:IE:uh} and~\eqref{eq:conv:IE:uhk}, respectively, we obtain
	\begin{equation}\label{eq:conv:IE:eh}
		\errh^{n+1} \le \beta\,C_F\,\errh^{n} + \sum_{k=1}^{m}\gamma_k\errhk^{n} + \O(m\dt(\dt + \dt^{\alpha}/m)),
	\end{equation}
	and
	\begin{equation}\label{eq:conv:IE:ehk}
		\errhk^{n+1} \le \gamma_k\errhk^{n} + \beta_k\,C_F\,\errh^{n+1} + \O(\dt^2),
		\qquad
		\errhk^{0} = 0.
	\end{equation}
	Recursive substitution in~\eqref{eq:conv:IE:ehk} yields to
	\begin{equation}\label{key}
		\errhk^{n} \le \beta_k\,C_F\,\sum_{j=1}^{n}\gamma_k^{n-j}\errh^{j} + \O(\dt).
	\end{equation}
	And substituting this to~\eqref{eq:conv:IE:eh}, we end up with
	\begin{equation}\label{key}
		\errh^{n+1} \le \beta\,C_F\errh^{n} + C_F\sum_{k=1}^{m}\sum_{j=1}^{n}\beta_k\gamma_k^{n+1-j}\errh^{j} + \O(m\dt),
		\qquad
		\errh^{0} = 0.
	\end{equation}
	Hence, the discrete Gr\"onwall inequality completes the proof:
	\begin{equation}\label{key}
		\errh^{n} \le \O(m\dt) \cdot \exp{C_F\left[\sum_{k=1}^{m}\sum_{j=0}^{n-1}\beta_k\gamma_k^{j} + w_{\infty}\right]} \le \O(m\dt) \cdot \exp{C_F \left[(t^n)^\alpha+\eps_{ra}\right]},
	\end{equation}
	where we used the following bound for the exponent:
	\begin{equation}\label{key}
		\sum_{j=0}^{n-1}\beta_k\gamma_k^{j}
		= \beta_k\frac{1-\gamma_k^n}{1 - \gamma_k}
		= w_k\frac{1-\gamma_k^n}{\lambda_k}
		\le \frac{w_k}{\lambda_k}\left(1 - \frac{1}{1 + n \lambda_k h}\right)
		\le \frac{w_k n\dt}{1+\lambda_k n\dt}.	
	\end{equation}
	$ $
\end{proof}

As we observed in~\Cref{fig:ModesNumber} that the number of modes $m=\O(\log\dt)$, the error of order $\O(mh)$ behaves asymptotically almost linearly.



\begin{remark}\label{rmk:CN:oscillations}
	Remark that the case $\theta=1/2$ in~\Cref{th:scheme}, the Crank-Nicolson scheme (CN), is known to be not L-stable.
	Indeed, we have $-1\le\gamma_k\le 1$ in~\eqref{eq:scheme:coefs} with $\theta=1/2$, moreover, $\gamma_k$ approaches $-1$ for large enough~$\lambda_k$, giving rise to a stiff problem.
	Thus, the higher modes produce undesired oscillation of the solution.
	When $\alpha$ decays, the tail of the function $z^{-\alpha}$ becomes heavier, and thus $\max \lambda_k$ grows (see, e.g.,~\Cref{fig:poles_bounds}).
	Therefore, the oscillations become more dominant, the smaller~$\alpha$ is.
	An example can be found in the next section (\Cref{fig:Stability}).
	This observation motivates us to introduce in the following lemma an alternative two-point scheme which expresses more stability but preserves the order.
\end{remark}




\begin{lemma}[Exponential Integrator]\label{th:exp:scheme}
	Let $\uht^{n}$ be defined by the following time-stepping scheme:
	\begin{equation}\label{eq:exp:final_scheme}
	\uht^{n+1} - \beta\,F_{-}[\uht^{n+1}] = \beta\, F_{+}[\uht^{n}] + u_0 + \sum_{k=1}^{m}\gamma_k\,\uhk^{n},
	\end{equation}
	where the discrete modes $u_k^{n}$, $k=1,\ldots,m$, are updated by
	\begin{equation}\label{eq:exp:final_scheme_modes}
	\begin{aligned}	
	\uhk^{n+1} = \gamma_k\,\uhk^{n} + \beta^1_k\,F[\tilde{u}_h^{n+1}] + \beta^2_k\, F[\uht^{n}],
	\qquad
	\uhk^{0} = 0,
	\end{aligned}
	\end{equation}
	with coefficients
	\begin{equation}\label{eq:exp:coefs}	
	\begin{gathered}
	\gamma_k = \exp{-\lambda_k\dt}, \qquad
	\beta^1_k = w_k\frac{\gamma_k - (1 - \lambda_k\dt)}{\lambda_k^2\dt}, \qquad
	\beta^2_k = w_k\frac{1 - (1 + \lambda_k\dt)\gamma_k}{\lambda_k^2\dt},
	\\	
	\beta = \sum_{k=1}^{m}(\beta^1_k + \beta^2_k) + w_{\infty} = \sum_{k=1}^{m}\frac{w_k}{\lambda_k}(1-\gamma_k) + w_{\infty}.
	\end{gathered}
	\end{equation}
	Then, $\uht^{n}$ approximates $\tilde{u}(t^n)$ with the discretization error of order~$\dt^{1+\alpha}$,
	\begin{equation}\label{key}
	\errh^{n} = \normH{\tilde{u}(t^{n}) - \uht^{n}} \le \O(h^{1+\alpha}) \cdot \exp{2\,C_F\left[(t^n)^\alpha + \eps_{ra}\right]}.
	\end{equation}
	where $\eps_{ra}:=\max\limits_{s\in[\frac{1}{T},\frac{1}{\dt}]}\left|\hat{K}(s)-\hat{\tilde{K}}(s)\right|$ stands for the rational approximation error.
\end{lemma}

\begin{proof}
	From~\eqref{eq:sys_modes}, the modes satisfy the recurrence relation
	\begin{equation}\label{eq:exp:uk_proof}
	u_k(t^{n+1}) = u_k(t^{n})\exp{-\lambda_k\dt} + \,w_k \int\limits_{t^{n}}^{t^{n+1}}\exp{-\lambda_k(t^{n+1}-s)}F[\tilde{u}](s)\d s.
	\end{equation}
	Let us consider the following quadrature rule for an arbitrary function~$f(t)$ and scalar~$\lambda\ge0$:
	\begin{align}
	\int\limits_{t^{n}}^{t^{n+1}}\exp{-\lambda(t^{n+1}-s)}f(s)\d s
	&= 
	a_1\frac{f(t^{n+1})+f(t^{n})}{2} 
	+ a_2 \frac{f(t^{n+1})-f(t^{n})}{\dt} 
	+ \O(a_1\dt^2)
	\\
	&=
	\left(\frac{a_1}{2}+\frac{a_2}{\dt}\right)f(t^{n+1})
	+ \left(\frac{a_1}{2}-\frac{a_2}{\dt}\right) f(t^{n})
	+ \O(a_1\dt^2)
	\end{align}
	with coefficients given as
	\begin{equation}\label{eq:Exp:int_coefs}
	a_1 = a_1(\lambda) = \int\limits_{0}^{\dt}\exp{-\lambda s}\d s \le \frac{2\dt}{1+\lambda\dt},
	\qquad
	a_2 = a_2(\lambda) = \int\limits_{0}^{\dt}\exp{-\lambda s}\left(\frac{\dt}{2}-s\right)\d s
	\end{equation}
	Applying the above quadrature rule to~\eqref{eq:exp:uk_proof}, given $\beta_k^1 := w_k\left(\frac{a_1(\lambda_k)}{2}+\frac{a_2(\lambda_k)}{\dt}\right)$ and $\beta_k^2 = w_k\left(\frac{a_1(\lambda_k)}{2}-\frac{a_2(\lambda_k)}{\dt}\right)$, we obtain
	\begin{equation}\label{eq:conv:Exp:uhk}
	u_k(t^{n+1}) = \gamma_k u_k(t^{n})
	+ \beta^1_k F[\tilde{u}](t^{n+1})
	+ \beta^2_k F[\tilde{u}](t^{n})
	+ \O(w_k a_1(\lambda_k)\dt^2).
	\end{equation}
	Note that $0\le a_1(\lambda_k) \le \frac{2\dt}{1+\lambda_k\dt}$ and thus $0 \le \sum_{k=1}^{m}w_k\,a_1(\lambda_k) \le 2 h^\alpha$.
	Summing up the modes, $u_0$ and $u_{\infty}$, we obtain
	\begin{equation}
		\tilde{u}(t^{n+1}) - \beta^1\,F[\tilde{u}](t^{n+1}) = \beta^2\,F[\tilde{u}](t^{n}) + u_0 + \sum_{k=1}^{m}\gamma_k u_k(t^{n}) + \O(\dt^{2+\alpha}),
	\end{equation}
	where 
	\begin{equation}\label{eq:beta12}
		\beta^1 := \sum_{k=1}^{m}\beta_k^1 +w_{\infty}
		\qquad\text{and}\qquad
		\beta^2 := \sum_{k=1}^{m}\beta_k^2.
	\end{equation}
	Note that $\beta_k^1 + \beta_k^2 \le 2\frac{w_k h}{1 + \lambda_k h}$ and thus $\beta^1 + \beta^1 = \O(h^\alpha)$.
	Then, using Taylor expansion of $F_{+}[\tilde{u}](t^{n+1})$ at the point $t^{n}$ and respectively $F_{-}[\tilde{u}](t^{n})$ at the point $t^{n+1}$, we write
	\begin{equation}\label{eq:conv:Exp:uh}
	\tilde{u}(t^{n+1}) - \beta\,F_{-}[\tilde{u}](t^{n+1}) = \beta\,F_{+}[\tilde{u}](t^{n}) + u_0 + \sum_{k=1}^{m}\gamma_k u_k(t^{n}) + \O(h^{1+\alpha}).
	\end{equation}
	We subtract~\eqref{eq:exp:final_scheme} and~\eqref{eq:exp:final_scheme_modes} from~\eqref{eq:conv:Exp:uh} and~\eqref{eq:conv:Exp:uhk}, respectively, to obtain
	\begin{equation}\label{eq:conv:Exp:eh}
	\errh^{n+1} \le \beta\,C_F\,\errh^{n} + \sum_{k=1}^{m}\gamma_k\errhk^{n} + \O(\dt^{1+\alpha}),
	\end{equation}
	and
	\begin{equation}\label{eq:conv:Exp:ehk}
	\errhk^{n+1} \le \gamma_k\errhk^{n} + C_F\,\left(\beta_k^1\,\errh^{n+1} + \beta_k^2\,\errh^{n}\right) + \O\left(\frac{w_k\dt^3}{1+\lambda_k h}\right),
	\qquad
	\errhk^{0} = 0.
	\end{equation}
	where we used continuity of~$F$ and monotonicity of~$F_{-}$ as in~\eqref{eq:proof:continuity_and_monotonicity}.
	Recursive substitution in~\eqref{eq:conv:Exp:ehk}, yields to
	\begin{equation}\label{key}
		\errhk^{n} \le C_F\,\sum_{j=1}^{n}(\beta_k^1+\beta_k^2)\gamma_k^{n-j}\errh^{j} + \O\left(\frac{w_k\dt^2}{1+\lambda_k h}\right).
	\end{equation}
	Substituting this to~\eqref{eq:conv:Exp:eh}, we thus write the estimation
	\begin{equation}\label{key}
	\errh^{n+1} \le C_F\left(\beta\,\errh^{n} + \sum_{k=1}^{m}\sum_{j=1}^{n}(\beta_k^1+\beta_k^2)\gamma_k^{n+1-j}\errh^{j}\right) + \O(h^{1+\alpha}),
	\qquad
	\errh^{0} = 0.
	\end{equation}
	Hence, by the discrete Gr\"onwall inequality, we finally obtain
	\begin{equation}\label{key}
		\errh^{n} 
		\le \O(h^{1+\alpha}) \cdot \exp{C_F\left[\sum_{k=1}^{m}\sum_{j=0}^{n-1}(\beta_k^1+\beta_k^2)\gamma_k^{j} + w_{\infty}\right]}
		\le \O(h^{1+\alpha}) \cdot \exp{2\,C_F\left[(t^n)^\alpha + \eps_{ra}\right]},
	\end{equation}
	using the following bound for the exponent:
	\begin{equation}\label{key}
	\sum_{j=0}^{n-1}(\beta_k^1+\beta_k^2)\gamma_k^{j}
	= \frac{w_k}{\lambda_k}(1-\gamma_k)\frac{1-\gamma_k^n}{1 - \gamma_k}
	= w_k \int\limits_{0}^{t^n}\exp{-\lambda_k t} \d t
	\le 2\,\frac{w_k t^n}{1+\lambda_k t^n}.	
	\end{equation}
	$ $
\end{proof}



\begin{remark}
	Note that the schemes presented in~\Cref{th:IE} and \Cref{th:exp:scheme} are of the implicit-explicit (IMEX) type.
	In particular, the explicit modal system defines the order of the scheme, while the implicit part, summing up the modes to the solution, guarantees stability via splitting of the operator.
\end{remark}

\begin{figure}[ht!]
	\centering\noindent
	\includegraphics[width=0.7\textwidth]{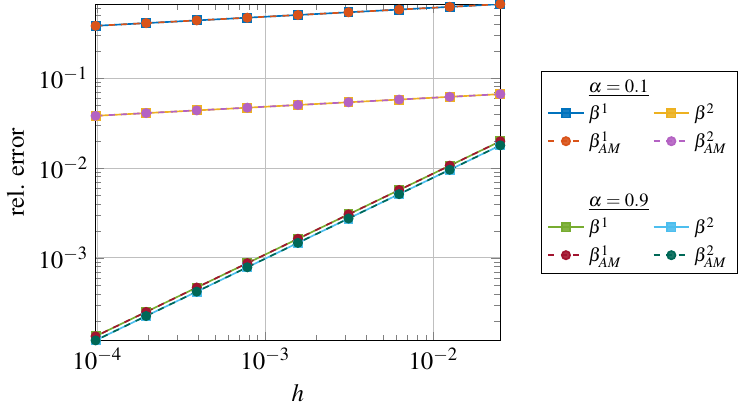}	
	\caption{Comparison of the coefficients $\beta^1$ and~$\beta^2$ in~\eqref{eq:exp:coefs} for the Exponential Integrator scheme (Lemma~\ref{th:exp:scheme}) with the $2$nd order fractional Adams--Moulton coefficients $\beta_{AM}^1=\dt^\alpha / \Gamma(\alpha+2)$ and $\beta_{AM}^2=\alpha \dt^\alpha/\Gamma(\alpha+2)$, respectively.}\label{fig:AM_coefs}
\end{figure}

\begin{remark}
	In contrast to the common strategy, when the fractional integral is split into the local and the history integrals (see, e.g., \cite{zayernouri2016fractional,baffet2017kernel,zhou2020implicit}), we did not discretize the local integral explicitly in the construction of our schemes.
	Instead, the local term is obtained as a linear combination of the modes (including the "infinity" mode~$u_{\infty}$ which does not however enter to the history part).
	Moreover, let us remark that the coefficients~$\beta^1$ and~$\beta^2$ in~\eqref{eq:beta12} approximate the second-order fractional Adams--Moulton coefficients~\cite{zayernouri2016fractional}, $\beta_{AM}^1=\dt^\alpha / \Gamma(\alpha+2)$ and $\beta_{AM}^2=\alpha \dt^\alpha/\Gamma(\alpha+2)$, respectively.
	The values of the coefficients for $\alpha=0.1$ and $\alpha=0.9$ are compared in~\Cref{fig:AM_coefs}.
	Thus, using Adams--Moulton type discretization for the modal equations, the rational approximation approach can automatically reconstruct the fractional Adams--Moulton coefficients, naturally leading to the local integration term arising in fractional linear multi-step methods.
\end{remark}

	\section{Numerical examples}
	\label{sec:Examples}


In this section, we illustrate the proposed scheme in application to the two following examples.
We first consider a simple linear case, more precisely, the one-dimensional fractional heat equation, where the analytical solution is known and given by a Mittag-Leffler function, so that we can study the accuracy and convergence rate.
Then, the scheme is applied to the more complex non-linear Cahn-Hilliard equation and compared to a classical fractional time-stepping scheme.
Both problems are discretized in space with Finite Elements using the~\fenics{} package~\cite{FEniCS}.
Since we focus on the accuracy of the time-integration scheme, we fix in what follows the space discretization to be sufficiently fine for not polluting the total error.	
	
		\subsection{Fractional heat equation}
		\label{sec:Linear}


Let us consider the one-dimensional fractional heat equation with homogeneous Dirichlet boundary conditions:
\begin{align}\label{eq:ScalarEq}
\drv{t}{^\alpha u}(t,x) - \drv{x}{^2 u}(t,x) &= 0, \qquad t\in(0,T], \; x\in(0, 1),\\
u(0,x) &= \sin(\pi\,x), \\
u(t,0) &= u(t,1) = 0.
\end{align}
Its analytical solution is given by $u(t,x) = \Ea_{\alpha,1}\left[-\pi^2\,t^\alpha\right]\,\sin(\pi\,x)$, see, e.g., \cite{Kexue2011}, where $\Ea_{\alpha,\beta}[x]$ is the Mittag-Leffler function~\eqref{eq:Mittag-Leffler}.
Note that in this example, we have $F[v] = F_{-}[v] = \drv{x}{^2 v}$, i.e. $F_{+}[v]\equiv 0$.
Thus, the scheme in \Cref{th:IE} coincides with the implicit Euler scheme in Proposition~\ref{th:scheme} with $\theta=1$.
For the spatial discretization, we use $5000$ $P_1$-elements.
We fix the final time $T=1$.
For the rational approximation, we set the AAA-tolerance~$10^{-13}$ with $100$ candidate points.
The dependence of the number of the modes~$m$ on the time step size~$\dt$ is shown in~\Cref{fig:ModesNumber}.

\begin{figure}[!ht]
	\centering\noindent
	\includegraphics[width=\textwidth]{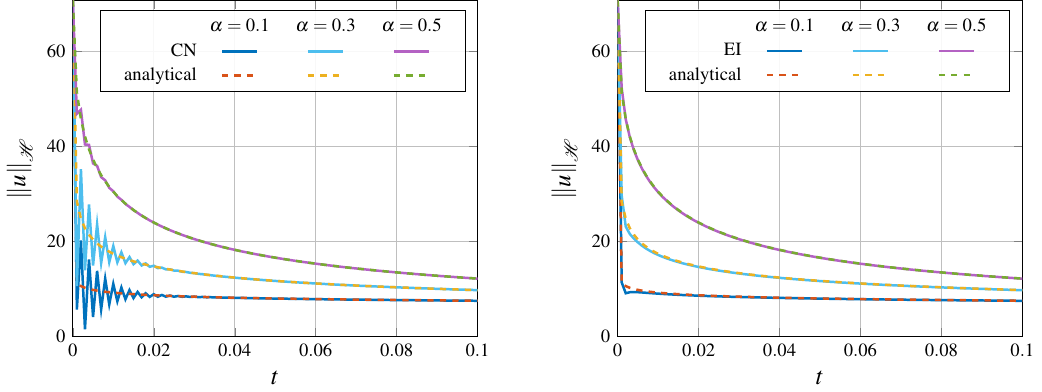}	
	
	\caption{Illustration of instability effects of the Crank-Nicolson (CN) scheme~(left) for $\alpha=0.1$, $0.3$, $0.5$ with $h=10^{-3}$.
		The effect is stronger for smaller~$\alpha$.
		However, there is no oscillatory effect in case of the Exponential Integrator (IE) scheme~\eqref{eq:exp:final_scheme}-\eqref{eq:exp:final_scheme_modes} (right). }
	\label{fig:Stability}
\end{figure}

We start with a comparison of the Crank-Nicolson (CN) scheme with the Exponential Integrator (EI) scheme (\Cref{th:exp:scheme}).
The evolution in time of the corresponding solutions with the time step size $\dt=10^{-3}$ for the cases~$\alpha=0.1$, $0.3$, $0.5$ is shown on~\Cref{fig:Stability}, where the norm $\normH{\cdot}$ stands for $\L{2}$-norm in space (using FE).
We observe oscillations for the CN scheme (left) but not for the EI scheme (right).
Moreover, the oscillations become stronger when $\alpha$ decreases.

\begin{figure}[ht!]
	\centering\noindent%
	\includegraphics[width=\textwidth]{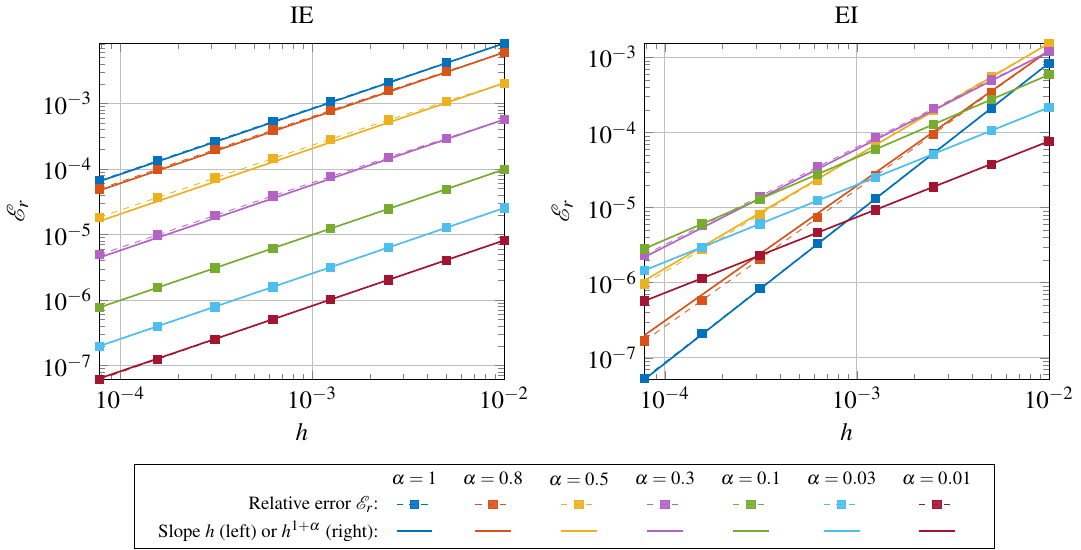}	
	\caption{Convergence rate of the Implicit Euler (IE) scheme~\eqref{eq:IE:scheme}-\eqref{eq:IE:modes} (left) and the Exponential Integrator (EI) scheme~\eqref{eq:exp:final_scheme}-\eqref{eq:exp:final_scheme_modes} (right). }
	\label{fig:ConvergenceRate}
\end{figure}

In view of \Cref{th:error_estimate}, we measure the global error using a weighted norm.
Thus, we consider the relative error $\ErrRel := {\norm{u-\uht}_{\alpha}/{\norm{u - u_0}_{\alpha}}}$,
where the norm $\norm{\cdot}_{\alpha}$ is defined via $\ell_2$-norm in time, weighted with $t^{2\cdot(1-\alpha)}$, and $\L{2}$-norm in space.	
We study the convergence rate of the error~$\ErrRel$ with respect to the time step size~$\dt$ for the Implicit Euler scheme \eqref{eq:IE:scheme}-\eqref{eq:IE:modes} and the Exponential Integrator scheme \eqref{eq:exp:final_scheme}-\eqref{eq:exp:final_scheme_modes}.
The results for $\alpha=1$, $0.8$, $0.5$, $0.3$, $0.1$, $0.03$, $0.01$ are plotted in~\Cref{fig:ConvergenceRate} and confirm the theoretical error bounds suggested in~\Cref{th:IE} and \Cref{th:exp:scheme}, respectively.
In particular, for the first scheme, we clearly observe linear convergence rate for all~$\alpha$, while the second scheme presents convergence of order $\dt^{1+\alpha}$.
The number of modes~$m$ varies between $8$ and $20$, according to~\Cref{fig:ModesNumber}, except for $\alpha=1$, when there is only one mode.

		\subsection{Fractional Cahn-Hillard equation}
		\label{sec:CahnHilliard}

Let $\Omega=(0,1)^2$ be a unit square domain and $\Bochner = \L{2}([0,T]; \HS)$, where $\HS=\H{1}(\Omega)\times\H{1}(\Omega)$.
Moreover, let~$\sca{\cdot}{\cdot}$ and $\norm{\cdot}$ denote the scalar product and the norm in $\L{2}(\Omega)$, respectively.
Then, we formulate the following non-linear Cahn-Hilliard problem:
find $(u, \mu) \in \Bochner$ satisfying for all $(v, w) \in \HS$ and all $t\in(0,T]$,
\begin{equation}\label{eq:CH}
\begin{aligned}
\sca{\drv{t}{^\alpha u}}{v} &= -\sca{M\grad\mu}{\grad v} 
,	\\
\sca{\mu}{w} &= \sca{\psi(u)}{w} + \eps^2\sca{\grad u}{\grad w},
\end{aligned}
\end{equation}
provided homogeneous Neumann boundary conditions, the initial state~$u(0)=u_0$, the constant mobility~$M$ and the surface parameter~$\eps$.
The non-linear function~$\psi(u)$ is defined as derivative of the potential $\Phi(u) = \frac{1}{4}(u^2-1)^2$:
\begin{equation}\label{key}
\psi(u) = \Phi^\prime(u) = u^3 - u.
\end{equation}
The Ginzburg--Landau free energy of the system is defined as
\begin{equation}\label{key}
E(u) = \int\limits_{\Omega}\left[\Phi(u) + \frac{\eps^2}{2}\abs{\grad u}^2\right] \d\Omega.
\end{equation}
Note that the function $\psi(u)$ is Lipschitz continuous on the interval~$[-1,1]$, i.e., between zeros of the potential~$\Phi(u)$.
Therefore, if the initial conditions are contained in the interval, the problem~$\eqref{eq:CH}$ satisfies the conditions of Theorem~\ref{th:error_estimate}.

Due to the "double-well" structure of the potential, presenting both convex and concave parts, stability of time-schemes for Cahn-Hilliard equation is a sophisticated question and is a subject of numerous works.
The fully implicit time-schemes are only conditionally stable~\cite{elliott1989cahn}.
Unconditionally stable schemes include so-called gradient stability, providing monotone decay of the discretized Ginzburg--Landau energy, e.g., splitting to implicit convex and explicit concave parts~\cite{eyre1998unconditionally,wu2014stabilized} or others~\cite{du1991numerical,gomez2011provably}.
The situation becomes more complicated in the case of fractional derivative, since even for the analytical solution, the associated Ginzburg--Landau energy is not proved to be monotone~\cite{tang2019energy}.
Note that the schemes presented in Lemmas~\ref{th:IE} and~\ref{th:exp:scheme} naturally allow splitting techniques, ensuring stability.
Splitting the potential into convex and concave parts, we consider $\psi(u) = \psi_{+}(u) + \psi_{-}(u)$, where the functions $\psi_{+}(u)$ and $\psi_{-}(u)$ are monotonously increasing and decreasing in~$[-1,1]$, respectively.
Note that such splitting is not unique.
Let us use the following splitting scheme, proposed in~\cite{eyre1998unconditionally2}:
\begin{gather}\label{eq:CH:split}
	\psi_{+}(u) = 2\,u,
	\qquad
	\psi_{-}(u) = u^3 - 3\,u.
\end{gather}

Let $\HS$ be an appropriate finite elements space.
Implementing the Exponential Integrator scheme~\eqref{eq:exp:final_scheme}-\eqref{eq:exp:coefs} for discretization of the problem~\eqref{eq:CH}, we compute at each time step~$n$ the discrete solution pair $(u_h^{n+1},\mu_h^{n+1}) \in \HS$ satisfying
\begin{align}\label{eq:CH:scheme}
\sca{u_h^{n+1}-H^n}{v} &= -\beta\,\sca{M\grad\tilde{\mu}_h^{n+1}}{\grad v} 
,	\\
\sca{\tilde{\mu}_h^{n+1}}{w} &= \sca{\psi_{+}(u_h^{n+1})+\psi_{-}(u_h^{n})}{w} + \eps^2\sca{\grad u_h^{n+1}}{\grad w},
\end{align}
for all $(v,w)\in\HS$, with the history term defined as $H^n = u_0 + \sum_{k=1}^{m}\gamma_k\,\uhk^{n}$, where the modes $u_k^{n}$, $k=1,\ldots,m$, are updated as follows:
\begin{equation}\label{eq:CH:scheme_modes}
\begin{aligned}	
\sca{\uhk^{n+1}}{\,v\,} &= \gamma_k\,\sca{\uhk^{n}}{v} -\sca{M\grad\left[\beta_k^1\mu^{n+1} + \beta_k^2\mu^{n}\right]}{\grad v},
\qquad
\uhk^{0} = 0,
\\
\sca{\mu_h^{n+1}}{w} &= \sca{\psi(u_h^{n+1})}{w} + \eps^2\sca{\grad u_h^{n+1}}{\grad w},
\end{aligned}
\end{equation}
with the coefficients $\beta$, $\beta_k^1$, $\beta_k^2$ and $\gamma_k$ given in~\eqref{eq:exp:coefs}.
We remark that a linear choice of the increasing part in~\eqref{eq:CH:split} leads to a linear implicit part in the time-scheme.
That is, though the problem~\eqref{eq:CH} is non-linear, the numerical solution of~\eqref{eq:CH:scheme} requires only a linear solver at each time step.

\begin{figure}[!ht]
\centering
\newcommand{\imwidth}{0.2\textwidth}
\setlength{\tabcolsep}{8pt}
\begin{tabular}{m{0.005\textwidth}*{4}{c}}
	& $\alpha=0.1$ & $\alpha=0.3$ & $\alpha=0.5$ & $\alpha=0.9$ \\[-8pt]
	&  &  &  &  \\
	\rotatebox{90}{$t=0$  }& 
	\makecell{\includegraphics[width=\imwidth]{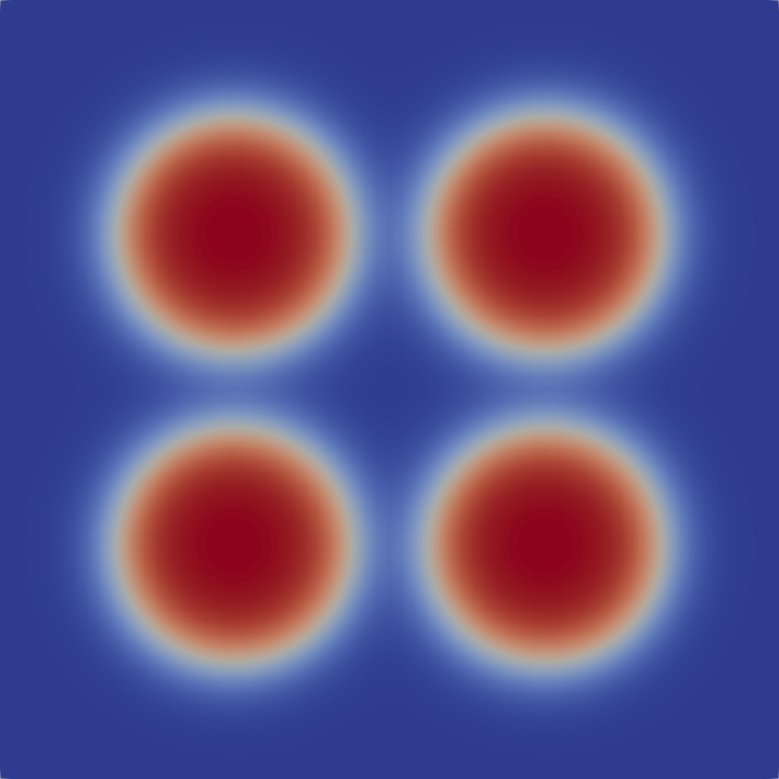}} &
	\makecell{\includegraphics[width=\imwidth]{figures/CahnHilliard/CH_a=0_1_t=0}} &
	\makecell{\includegraphics[width=\imwidth]{figures/CahnHilliard/CH_a=0_1_t=0}} &
	\makecell{\includegraphics[width=\imwidth]{figures/CahnHilliard/CH_a=0_1_t=0}} \\
	\rotatebox{90}{$t=0.4$  }& 
	\makecell{\includegraphics[width=\imwidth]{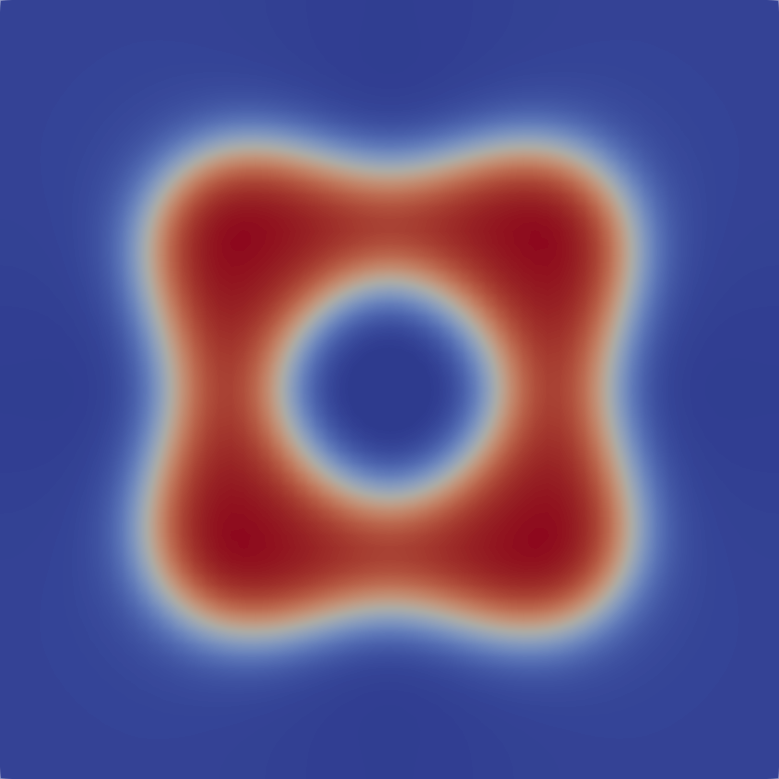}} &
	\makecell{\includegraphics[width=\imwidth]{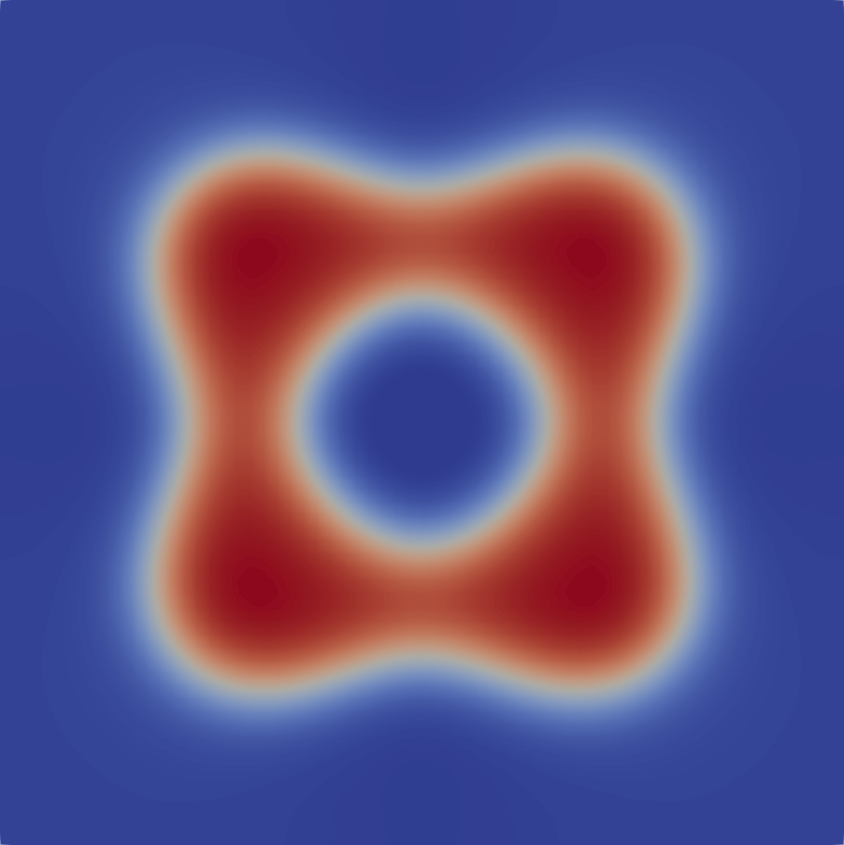}} &
	\makecell{\includegraphics[width=\imwidth]{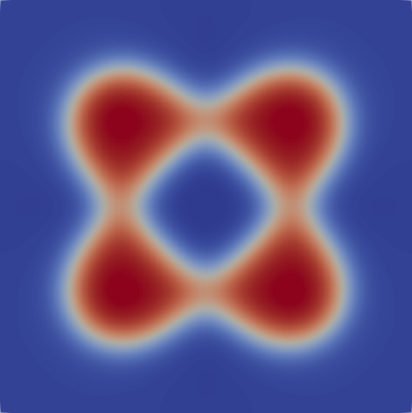}} &
	\makecell{\includegraphics[width=\imwidth]{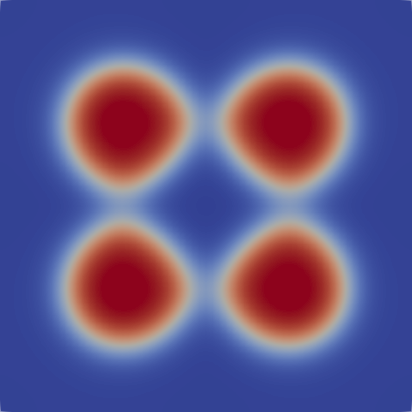}} \\
	\rotatebox{90}{$t=3.2$  }& 
	\makecell{\includegraphics[width=\imwidth]{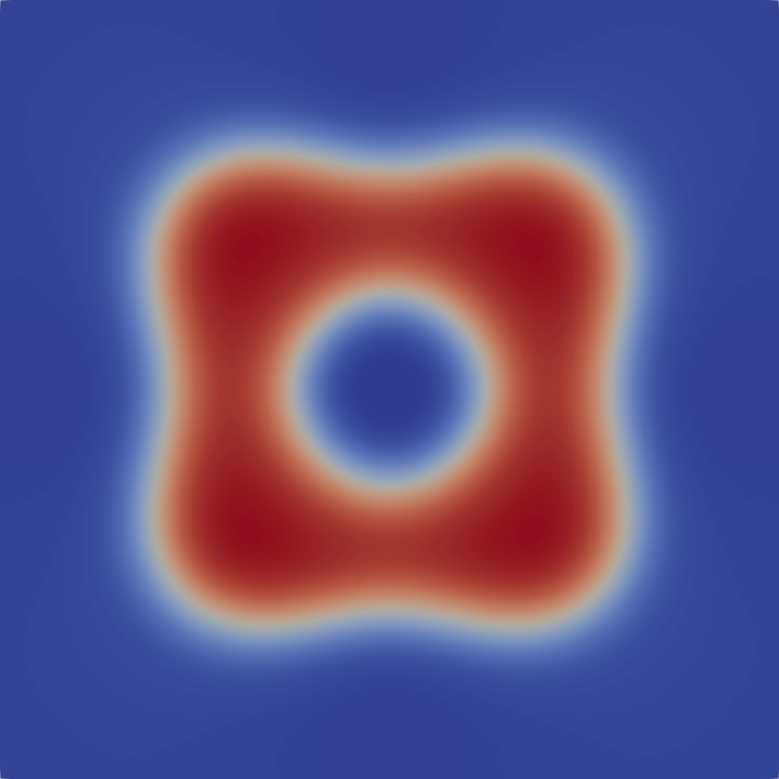}} &
	\makecell{\includegraphics[width=\imwidth]{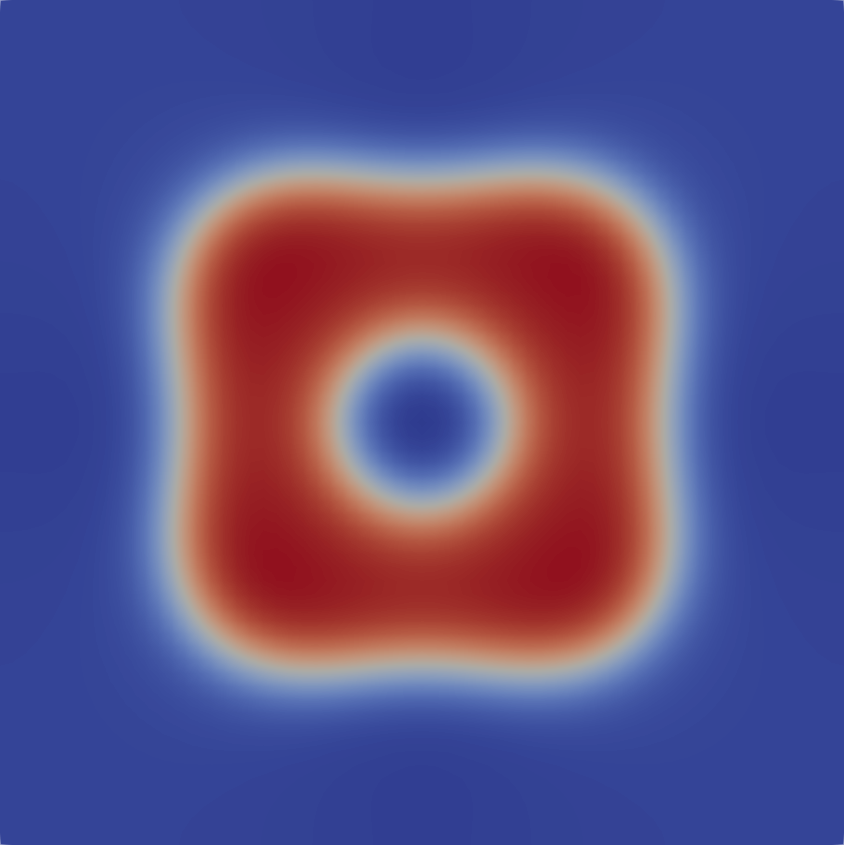}} &
	\makecell{\includegraphics[width=\imwidth]{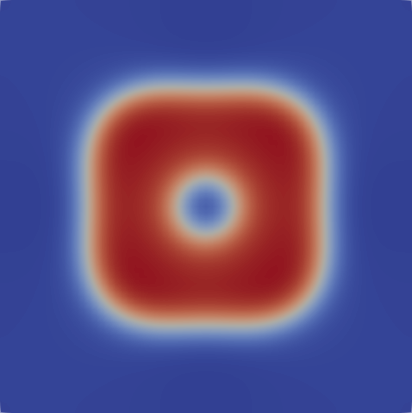}} &
	\makecell{\includegraphics[width=\imwidth]{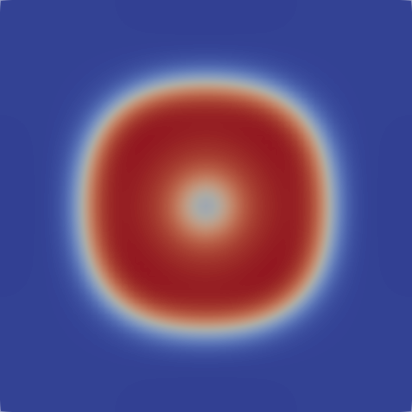}} \\
	\rotatebox{90}{$t=4$  }& 
	\makecell{\includegraphics[width=\imwidth]{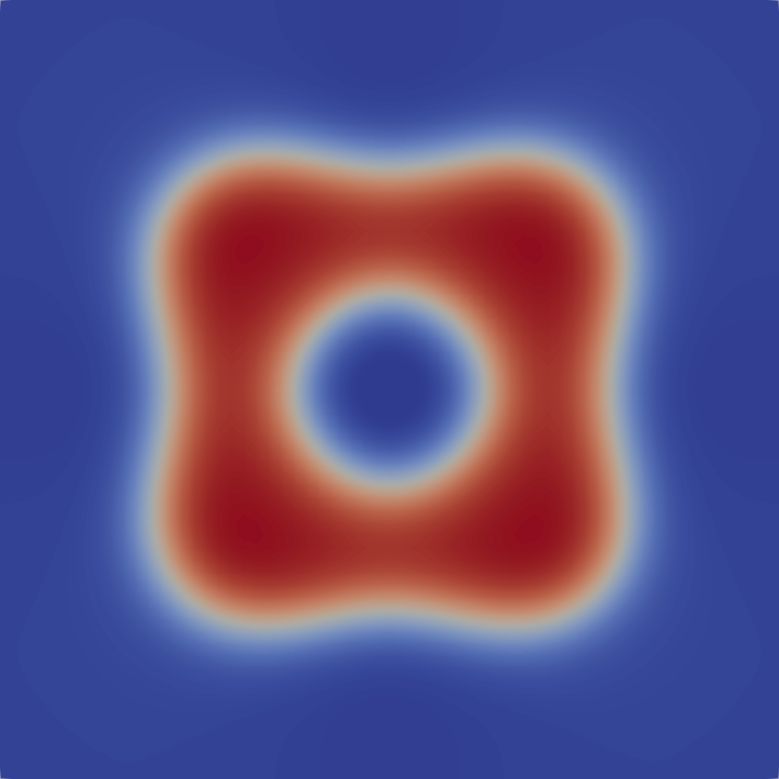}} &
	\makecell{\includegraphics[width=\imwidth]{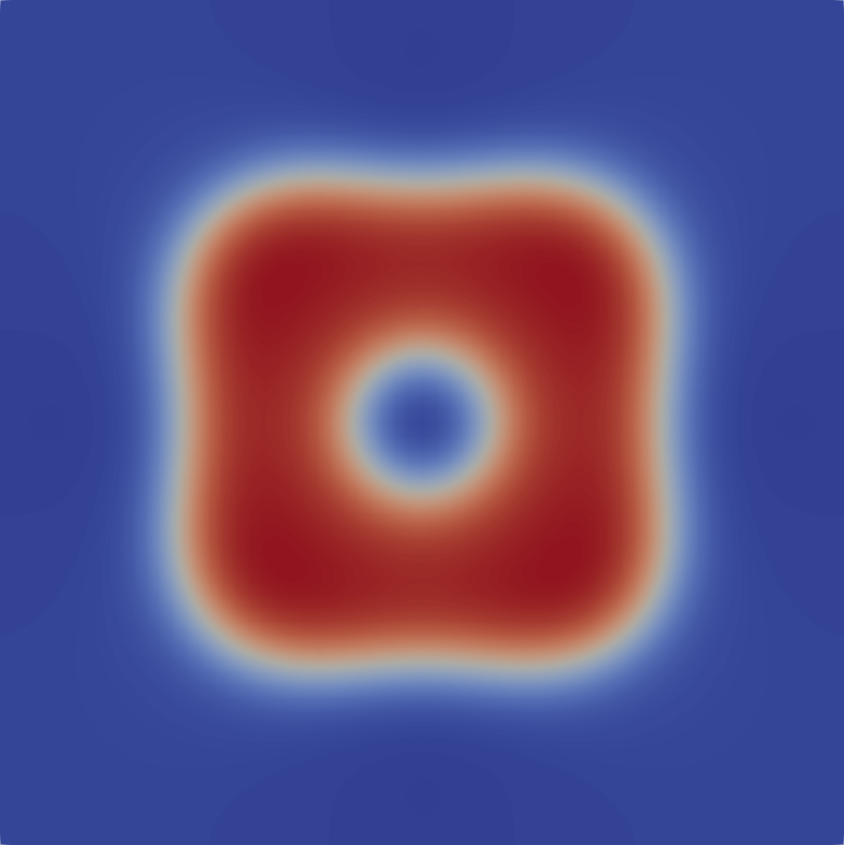}} &
	\makecell{\includegraphics[width=\imwidth]{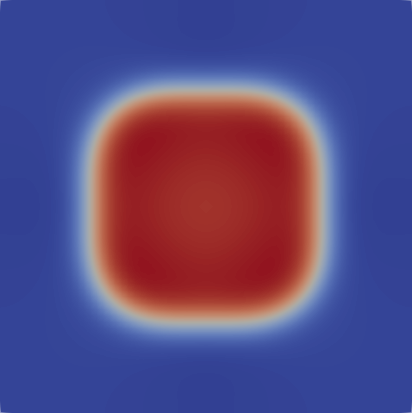}} &
	\makecell{\includegraphics[width=\imwidth]{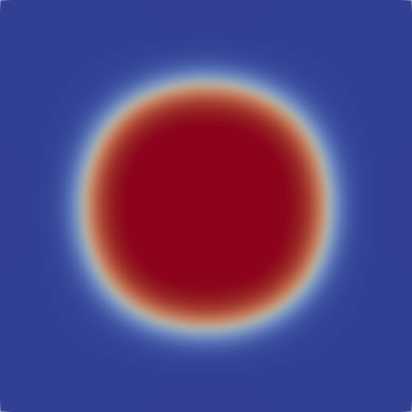}} \\
\end{tabular}
\caption{Fractional non-linear Cahn-Hilliard equation with $\alpha=0.1$, $0.3$, $0.5$, $0.9$.
	Four bubbles coalesce with different speed: for smaller $\alpha$, the coalescence accelerates in the beginning but then slows down with respect to larger values of $\alpha$.}
\label{fig:CH:Coalescence}
\end{figure}

\begin{figure}[!ht]
	\centering\noindent%
	\includegraphics[width=0.45\textwidth]{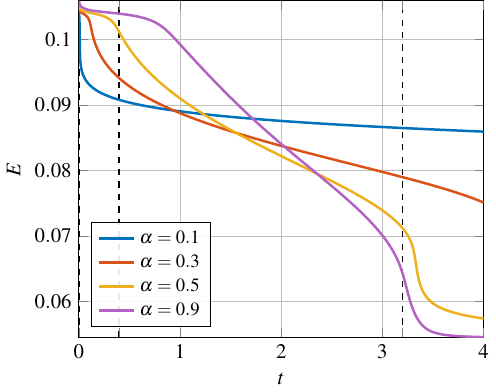}	
	\hspace{3ex}
	\centering\noindent%
	\includegraphics[width=0.45\textwidth]{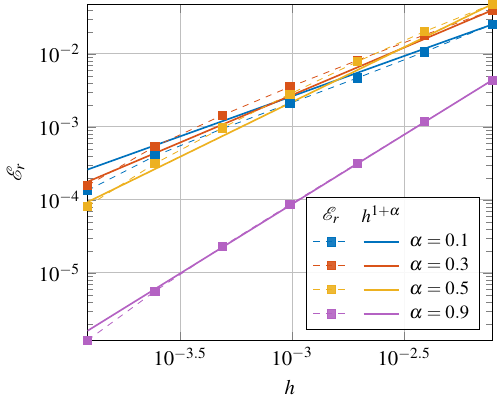}	
	\caption{Fractional non-linear Cahn-Hilliard equation with $\alpha=0.1$, $0.3$, $0.5$, $0.9$.
		\textbf{Left:} Evolution in time of the Ginzburg--Landau energy.
		Vertical dashed lines indicate the time points of the solution snapshots from~\Cref{fig:CH:Coalescence}.
		\textbf{Right:} Convergence rate for the rational approximation based scheme~\eqref{eq:CH:scheme}-\eqref{eq:CH:scheme_modes}.
	}
	\label{fig:CH:Energy_Convergence}
\end{figure}

For our simulation, we fix the constant mobility~$M=0.05$,  the surface parameter~$\eps=0.03$ and the final time~$T=4$.
For discretization in space, we use $(Q_1, Q_1)$ elements on a $64\times 64$ quadrilateral mesh.
For the rational approximation, we use the AAA-tolerance~$10^{-12}$ with $100$ candidate points (a logarithmic grid).
The initial state is
\begin{equation}\label{key}
u_0(\x) = \sum_{i=1}^{4}\tanh\frac{r-\abs{\x-\x_i}}{\sqrt{2}\eps} + 3,
\qquad \x\in\Omega,
\end{equation}
with $\x_1=(0.3,0.3)$, $\x_2=(0.3,0.7)$, $\x_3=(0.7,0.7)$, $\x_4=(0.7,0.3)$ and $r=0.15$, which corresponds to four bubbles of radius~$r$ centered at $\x_i$.
Due to the surface tension, the bubbles tend to coalesce in time~\cite{liu2018time}.
However, the process proceeds with different speeds for different values of the fractional order~$\alpha$.
In particular, for a small $\alpha$, the coalescence accelerates in the beginning but then slows down with respect to larger values of $\alpha$.
This effect can be observed in Figure~\ref{fig:CH:Coalescence}, where different states (computed with $\dt=2^{-14}$) are shown for $\alpha=0.1$, $0.3$, $0.5$, $0.9$ at time $t=0$, $0.4$, $3.2$, $4$.
Such behavior is also observed for the evolution of the corresponding Ginzburg--Landau energies, which is depicted in Figure~\ref{fig:CH:Energy_Convergence} on the left,
where the vertical dashed lines indicate the time points $t=0$, $0.4$, $3.2$, $4$ of the solution snapshots from~\Cref{fig:CH:Coalescence}.
These solutions are taken as reference for the convergence study of the relative error~$\ErrRel$.
In Figure~\ref{fig:CH:Energy_Convergence} on the right, there are plotted the convergence rates of the error with respect to the time step size~$h$ for the same values of $\alpha$.
We can observe that the error convergence respects the theoretical bounds.


	\section{Conclusion}
	\label{sec:Conclusion}

In this work, we proposed a new numerical method for solving fractional in time differential equation.
The method is based on the approximation of the Laplace spectrum of the fractional convolution kernel with a rational function, more precisely, a multi-pole series with an additional constant term.
To this end, we used the barycentric rational interpolation with the adaptive Antoulas--Anderson (AAA) algorithm.
This leads to the approximation of the kernel itself with a sum-of-exponentials with an additional singular term.
Thus, the solution of the FODE is represented as a sum of a small number of modes $m$ which solve a system of ODEs and can be updated in parallel.
The number of modes $m$ grows as $\log \tfrac{1}{\dt}$, leading to the complexity of order $\O(\tfrac{1}{\dt}\,\log \tfrac{1}{\dt})$ and memory requirements $\O(\log \tfrac{1}{\dt})$ with the time steps size~$\dt$, which is typical for numerical methods for time-fractional differential equations.
However, in our method, the value of the number of modes $m$ is significantly less than in many other kernel compression methods.
We proposed two new numerical time-integration schemes with convergence orders $\O(\dt\log\dt)$ and $\O(h^{1+\alpha})$.
The accuracy of the schemes is illustrated through the solution of a linear problem with known analytical solution.
The method is also applied to a non-linear fractional Cahn-Hilliard problem in 2D.

	
	\section*{Acknowledgements}
	
	This work was funded by the German Research Foundation by grants WO671/11-1 and the European Union's Horizon 2020 research and innovation programme under grant agreement No 800898.

	
	\bibliographystyle{IMANUM-BIB}
	\bibliography{bibliography}

\end{document}